\documentclass[11pt,a4paper]{amsart}
\usepackage[T1]{fontenc}
\usepackage[utf8]{inputenc}
\usepackage{lmodern,microtype,amsmath,amssymb,amsthm,mathtools}
\usepackage[a4paper,left=23mm,right=23mm,top=23mm,bottom=24mm]{geometry}
\usepackage{enumitem}
\usepackage{needspace}
\usepackage[hidelinks,breaklinks=true]{hyperref}
\setlist[itemize]{leftmargin=17pt,itemsep=3pt,topsep=3pt}
\newcommand{\R}{\mathbb R}\newcommand{\C}{\mathbb C}
\newcommand{\dd}{\mathrm d}\newcommand{\Dom}{\operatorname{dom}}\newcommand{\ran}{\operatorname{ran}}
\newcommand{\supp}{\operatorname{supp}}\newcommand{\rank}{\operatorname{rank}}
\newcommand{\diag}{\operatorname{diag}}\newcommand{\tr}{\operatorname{tr}}
\newcommand{\norm}[1]{\left\|#1\right\|}
\newcommand{\ip}[2]{\langle #1,#2\rangle}
\newtheorem{theorem}{Theorem}[section]
\newtheorem{proposition}[theorem]{Proposition}
\newtheorem{lemma}[theorem]{Lemma}
\newtheorem{corollary}[theorem]{Corollary}
\theoremstyle{remark}\newtheorem{remark}[theorem]{Remark}
\theoremstyle{definition}\newtheorem{problem}[theorem]{Open problem}

\hypersetup{pdftitle={Scalar realizability of biharmonic point interactions: critical scales and geometric obstructions},pdfauthor={Niyaz Tokmagambetov}}
\title[Scalar biharmonic realizability]{Scalar realizability of biharmonic point interactions: critical scales and geometric obstructions}
\author[N. Tokmagambetov]{Niyaz Tokmagambetov}
\address{Institute of Mathematics and Mathematical Modeling, Almaty, Kazakhstan}
\email{tokmagambetov@math.kz}
\email{niyaz.tokmagambetov@gmail.com}
\subjclass[2020]{Primary 47A55, 35J40; Secondary 47B25, 35P15}
\keywords{Biharmonic operator, point interaction, norm resolvent convergence, collision asymptotics, singular potential}
\numberwithin{equation}{section}
\allowdisplaybreaks
\begin{document}

\begin{abstract}
We study which planar biharmonic point interactions arise from colliding
scalar sites and fixed-profile smooth multiplication potentials with a
common lower spectral bound. We classify the two-site limits and derive
their critical smoothing laws and leading errors. For a symmetric
Hessian triplet we classify every separated smoothing width, including
oscillating widths, and identify two sharp logarithmic thresholds.
For the point models and their norm resolvent approximations by smooth
potentials, derivative closure and matrix inertia show that activating
both traceless Hessian channels requires at least five sites, at least
two of them attractive. A two-scale construction realizes both channels
with nine sites, three of them attractive. A full six-channel correction
requires at least seven sites and at least five attractive sites.
An obtuse triangle realizes every real symmetric dipole block,
together with a pinned monopole, with the minimal counts of three sites
and two attractive sites.
The operator framework is included: we specify the complete admissible
jets, repair a punctured ambient domain and calculate the full planar
Weyl matrix. 
\end{abstract}
\maketitle

\section{Introduction}
We consider the scalar forms
\[
q_h[u]=\|\Delta u\|_2^2+\sum_j\alpha_j(h)|u(ha_j)|^2,
\qquad u\in H^2(\R^2),\quad \alpha_j(h)\in\R,
\]
and their regularization by smooth multiplication potentials with one
fixed nonnegative profile. A lower spectral bound independent of $h$
is part of the approximation problem. Our question is which derivative
channels and coupling relations survive, and how the answer depends on
site geometry, signs, collision scales and smoothing width.

\paragraph{Extension theory and scalar realizability.}
Karpeshina and Pavlov classified the self-adjoint point extensions of
the planar biharmonic operator and identified the admissible multipoles
for the polyharmonic operator in general dimension
\cite{KP}. For a smooth submanifold of dimension $p$ in
$\R^n$, Fraguela Collar studied the capacity and trace conditions for
delta interactions and their strong resolvent approximation
\cite{Fra81}. For nonnegative smooth approximants to surface measure
with shrinking support, and $m>(n-p)/2$, both choices of sign have
a common lower bound and converge in strong resolvent sense.
The convergence is in norm when $m>n/2$
\cite{Fra87}. A single point with $m=n=2$ gives
the fixed-sign planar biharmonic monopole. Albeverio--Kurasov \cite{AK99,AK00} and
Posilicano \cite{Pos01,Pos} give the singular perturbation and
finite-trace calculus used here. Kurasov--Boman \cite{KB98} treat
one-dimensional derivative interactions in a different distributional
setting. Existence of a point extension does not prescribe diagonal
scalar approximants with a common lower bound.

The ambient-domain issue in \cite{KT} is
addressed in Appendix~\ref{fd:sec:found-restrictions}. The regularized
trace formulas of \cite{KT15} concern a related class of well-posed
perturbations and are not examined here. On the centered
Navier disk, $(r^4-4r^2+3)\cos2\theta$ is a nonzero bounded homogeneous
solution with all the listed radial and first angular fluxes zero.
Proposition~\ref{fd:prop:counterexample} identifies its Hessian source
and the exact two-dimensional invisible $L^2$ kernel in the declared
outer regularity class. Theorem~\ref{fd:thm:repair} specifies the
replacement direct sum and its exact $L^2$ criterion.
Theorem~\ref{fd:thm:planarM} evaluates the full Weyl matrix in one
fixed jet basis and anchor. We claim neither a new six-channel
classification nor priority for a matrix in an equivalent chart.
These appendices supply the domain correction and coordinates needed
to formulate the scalar obstruction results.

\paragraph{Approximation precedents.}
For one-dimensional Schr\"odinger operators, Cheon and Shigehara used
three approaching delta sites with separation-dependent strengths to
produce a jump in the wave function \cite{CS98}.
Exner, Neidhardt and Zagrebnov obtained a $\delta'$ interaction as a
norm resolvent limit using a tuned three-site configuration and scaled
short-range potentials \cite{ENZ}. Golovaty--Hryniv \cite{GH13} and Golovaty
\cite{Gol13} establish other one-dimensional norm resolvent limits,
including dependence on profile, scale ratio and resonance.
For a fourth-order operator on an interval, Man'ko showed how smooth
regularizations of a delta distribution and its derivatives select
different limiting operators \cite{Manko10}.

For the planar Laplacian, Fassari, Popov and Rinaldi proved norm
resolvent convergence of two coalescing, equally coupled point
interactions to a single attractive point interaction under a
distance-dependent renormalization \cite{FPR20}. For a Kirchhoff plate,
Haslinger, Jones, Movchan and Movchan used a dipole approximation to
describe scattering by closely spaced pairs of rigid pins \cite{HJMM18}.
Griesemer--Hofacker \cite{GH22} and Michelangeli--Scandone \cite{MS19}
study different contact and fractional-Laplacian approximation classes.
Here the sites enter a scalar biharmonic form, and we determine the
smoothing thresholds and site counts described below.

\paragraph{Main results.}
\begin{enumerate}
\item Theorems~\ref{thm:critical} and \ref{thm:converse} classify the
two-site limits in a fixed direction. The full-rank chart forces a
vanishing monopole diagonal parameter, and a rank-one dipole is
impossible. Theorems~\ref{thm:stable}, \ref{thm:smooth} and
\ref{thm:rate} give a common lower bound, smooth realizations and the
nonzero leading error.
\item Theorems~\ref{thm:pairwidth} and \ref{thm:pairtwoscale} give the
critical pair width and competing error orders. For a nonzero pair
parameter, $\varepsilon=h^p$ preserves the target for $p>1$, shifts
the coupling at $p=1$ and produces a point pin for $0<p<1$.
Theorem~\ref{thm:triplet} treats the symmetric Hessian triplet, its
common lower bound and logarithmic point-model error.
Theorem~\ref{thm:hessianwidth} covers every separated width,
including oscillating widths, and Corollary~\ref{cor:widththresholds}
gives two distinct sharp logarithmic thresholds. Table~\ref{tab:smoothing}
compares the thresholds for the pair and the triplet. The fixed smooth
profile may be asymmetric and need not be centred.
\item Theorem~\ref{thm:twohessians} realizes both traceless Hessians
with nine sites, three attractive sites and two scales. Exact angular
cancellation and finite-matrix inertia give the common lower bound
before either limit is taken.
\item Theorems~\ref{thm:derivativeclosure},
\ref{thm:sixsiteobstruction} and \ref{thm:attractiverank} give
geometric and sign obstructions: two Hessians require at least five
sites; a rank-six correction requires at least seven sites and five
attractive sites. Theorem~\ref{thm:obtusedipoles} realizes every real
symmetric dipole block, together with a pinned monopole, with the
minimal counts of three sites and two attractive sites.
\end{enumerate}
The exact Hessian minima and full real six-channel scalar reachability
are open. Necessary counts are not sufficiency theorems.

\begin{table}[!ht]
\caption{Sharp smoothing thresholds for the pair \eqref{eq:prescription}
and the symmetric Hessian triplet \eqref{eq:triscaling},
with $L=\log(1/h)$. Each column gives the condition for retaining
the indicated quantity of the point model.}
\label{tab:smoothing}
\centering
\small
\renewcommand{\arraystretch}{1.25}
\begin{tabular}{@{}lcc@{}}
\hline
Configuration & Limiting operator & Leading resolvent error term \\
\hline
Pair ($a\ne0$, $\varepsilon=h^p$, $p>0$) & $p>1$ & $p>3/2$ \\
Hessian triplet ($\varepsilon/h\to0$)
  & $\varepsilon=o(h^2/\sqrt L)$ & $\varepsilon=o(h^2/L)$ \\
\hline
\end{tabular}
\end{table}

\paragraph{Scope and organization.}
The main proofs concern whole-space scalar collisions with a declared
common lower bound. The bounded-domain framework in the appendices
does not by itself prove collision convergence for another outer
realization. Section~\ref{sec:loc-pairs} treats pairs;
Section~\ref{sec:loc-smoothing} treats smoothing and rates;
Section~\ref{sec:loc-hessian} treats Hessian widths;
Sections~\ref{sec:loc-obstructions}--\ref{sec:loc-multipoles}
give stability, geometry and minimum counts. The appendices provide the operator framework.

\section{Two-site point interactions}\label{sec:loc-pairs}

\subsection{Model and normalization}\label{sec:loc-model}

Fix a unit vector $e\in\R^2$ and put $H_0=\Delta^2$,
$\Dom(H_0)=H^4(\R^2)$. For finite nonzero real strengths $\alpha_j(h)$,
let $H_h$ be the form realization
\begin{equation}\label{eq:form}
 q_h[u]=\norm{\Delta u}_2^2+\alpha_1(h)|u(0)|^2
                         +\alpha_2(h)|u(he)|^2,
 \qquad \Dom(q_h)=H^2(\R^2).
\end{equation}
Write $b_j(h)=\alpha_j(h)^{-1}$. A zero strength can instead be treated
by deleting its site. Inner products are linear in the second argument.
We use
\[
 \widehat f(\xi)=\int_{\R^2}e^{-ix\cdot\xi}f(x)\,\dd x,
 \qquad f(x)=(2\pi)^{-2}\int_{\R^2}e^{ix\cdot\xi}\widehat f(\xi)\,\dd\xi.
\]
For $\lambda\in\C\setminus(-\infty,0]$, with principal roots and logarithm,
put
\begin{equation}\label{eq:green}
 G_\lambda(x)=(2\pi)^{-2}\int_{\R^2}
       \frac{e^{ix\cdot\xi}}{|\xi|^4+\lambda}\,\dd\xi,
 \qquad g_\lambda=G_\lambda(0)=\frac1{8\sqrt\lambda}.
\end{equation}
The value of $g_\lambda$ follows by radial integration at $\lambda>0$
and holomorphy.
Set $a_{0,t}[u]=\norm{\Delta u}_2^2+t\norm{u}_2^2$ for $t>0$.
Fourier Cauchy--Schwarz gives the uniform trace inequality
\begin{equation}\label{eq:trace}
 |u(x)|^2\le g_t a_{0,t}[u].
\end{equation}
Since $g_t\to0$ as $t\to\infty$, for fixed $h$ the point terms are form
bounded with respect to $\Delta^2$ with relative bound zero, and the KLMN
theorem \cite[Section~6.5]{Teschl} defines $H_h$ as a self-adjoint operator
bounded below. No lower bound uniform in $h$ is asserted here. The same
applies to any finite number of sites.

For sites $x_j$, let $B=\diag(b_j)$, $M_\lambda=(G_\lambda(x_i-x_j))_{i,j}$,
and $\Phi_\lambda c=\sum_jc_jG_\lambda(\cdot-x_j)$. Their resolvent is
\begin{equation}\label{eq:res}
 R_h(\lambda):=(H_h+\lambda)^{-1}
 =R_0(\lambda)-\Phi_\lambda(B+M_\lambda)^{-1}\Phi_{\bar\lambda}^*,
 \qquad R_0(\lambda)=(H_0+\lambda)^{-1}.
\end{equation}
Here $M_\lambda$ is the matrix of Green kernel values at the current
sites. At a positive real shift it is the Gram matrix of the evaluation
representers in $(H^2,a_{0,\lambda})$. In the extension formulas of
Appendices~\ref{fd:sec:found-extensions}--\ref{fd:sec:found-repairmatrix},
$M_z$ denotes the Weyl matrix for the chosen trace map. The spectral
parameters are related by $z=-\lambda$, so the planar free resolvents
satisfy $\mathcal R_0(-\lambda)=R_0(\lambda)$.

Formula \eqref{eq:res} specializes the finite-rank resolvent formula
\cite[Theorem~2.1]{Pos01} to point evaluations with inverse strengths $B$
(see also \cite[Chapter~3]{AK00}). Applying $H_0+\lambda$ to the
columns of $\Phi_\lambda$ gives the independent point masses. Thus
$M_\lambda>0$ for $\lambda>0$, and every nonzero $c\in\ker(B+M_\lambda)$
gives an eigenfunction $\Phi_\lambda c$ at $-\lambda$.

For the pair specifically,
\begin{equation}\label{eq:Q}
 \Phi_h(\lambda)=(G_\lambda,G_\lambda(\cdot-he)),\qquad
 Q_h(\lambda)=\begin{pmatrix}b_1+g_\lambda&G_\lambda(he)\\
 G_\lambda(he)&b_2+g_\lambda\end{pmatrix}.
\end{equation}
Real positive shifts in resolvent formulas are always restricted to their
resolvent sets. For self-adjoint operators it suffices to prove norm
convergence to the limiting resolvent at one common nonreal parameter
\cite[Section~6.6]{Teschl}. The resolvent identity gives the same
implication at a shared real resolvent point.

\begin{lemma}[Kernel estimates and renormalized moments]\label{lem:kernel}
Put $D_\lambda=-\partial_eG_\lambda$, $E_\lambda=\partial_e^2G_\lambda$,
and $D_{h,\lambda}=(G_\lambda(\cdot-he)-G_\lambda)/h$. Then
\begin{gather}
 G_\lambda\in H^2,\qquad D_\lambda\in H^1\setminus H^2,
 \qquad E_\lambda\in L^2\setminus H^1,\label{eq:sob}\\
 \norm{D_{h,\lambda}-D_\lambda}_2
 \le\tfrac h2\norm{\partial_e^2G_\lambda}_2.\label{eq:D}
\end{gather}
With $\ell_\lambda=\gamma_{\rm Eul}-\log2+\tfrac14\log\lambda$, where $\gamma_{\rm Eul}$ is
Euler's constant, the local expansion is
\begin{align}
 G_\lambda(re)=g_\lambda
 &+\frac{r^2}{8\pi}(\log r+\ell_\lambda-1)
       -\frac{\sqrt\lambda}{512}r^4\nonumber\\
 &+\frac{\lambda r^6}{4608\pi}
       (\tfrac{11}{6}-\ell_\lambda-\log r)+O(r^8).
                                                        \label{eq:local}
\end{align}
In particular,
\begin{align}
 \frac{2[(g_\lambda-G_\lambda(he))-(g_1-G_1(he))]}{h^2}
 &=-\frac{\log\lambda}{16\pi}
   +\frac{\sqrt\lambda-1}{256}h^2
   +O(h^4(1+|\log h|)).\label{eq:renorm}
\end{align}
All constants may depend on a fixed $\lambda$ in the cut plane.
\end{lemma}
\begin{proof}
The multiplier of an order-$j$ derivative of $G_\lambda$ has radial
order $j-4$ and a nonzero angular factor. Polar integration gives
\eqref{eq:sob}, with logarithmic divergence at each excluded endpoint.
Taylor's formula for translations in $L^2$ gives \eqref{eq:D}.
Put $F_\lambda=(|\xi|^4+\lambda)^{-1}-(|\xi|^4+1)^{-1}$. The second
and fourth moments needed in Appendix~\ref{fd:sec:found-weyl} are
\begin{align}
 I_2(\lambda)&=(2\pi)^{-2}\int_{\R^2}\xi_e^2F_\lambda(\xi)\,\dd\xi
 =\frac1{4\pi}\int_0^\infty r^3
 \left(\frac1{r^4+\lambda}-\frac1{r^4+1}\right)\dd r
 =-\frac{\log\lambda}{16\pi},\label{eq:I2}\\
 I_4(\lambda)&=(2\pi)^{-2}\int_{\R^2}\xi_e^4F_\lambda(\xi)\,\dd\xi
 =\frac3{16\pi}\int_0^\infty r^5
 \left(\frac1{r^4+\lambda}-\frac1{r^4+1}\right)\dd r
 =\frac3{64}(1-\sqrt\lambda),\label{eq:I4}
\end{align}
by the substitutions $t=r^4$ and $t=r^2$, respectively, followed by
analytic continuation from $\lambda>0$. For the full kernel, partial
fractions and the integral representation of $K_0$
\cite[Sections~9.5--9.6]{Temme} give
\begin{equation}\label{eq:K0complex}
 G_\lambda(re)=\frac{K_0(e^{-i\pi/4}\lambda^{1/4}r)
                         -K_0(e^{i\pi/4}\lambda^{1/4}r)}
                        {4\pi i\sqrt\lambda}.
\end{equation}
Both arguments have positive real part. In the convergent $K_0$ series
\cite[formula~(9.30)]{Temme}, terms of degrees $0,2,4,6$ give
\eqref{eq:local}, with harmonic numbers $1,3/2,11/6$. The logarithmic
terms of degree eight cancel between the two arguments, leaving an
$O(r^8)$ remainder. Subtracting the expansion at $\lambda=1$ gives
\eqref{eq:renorm}.
\end{proof}

\begin{lemma}[The full admissible point space]\label{lem:sixjets}
Any self-adjoint norm resolvent limit of scalar sites converging to zero
extends $\Delta^2|_{C_c^\infty(\R^2\setminus\{0\})}$.
At a nonreal shift its resolvent correction has range in
\[
 \operatorname{span}\{G_\lambda,\partial_1G_\lambda,\partial_2G_\lambda,
 \partial_1^2G_\lambda,\partial_1\partial_2G_\lambda,
                                      \partial_2^2G_\lambda\}.
\]
These six functions are independent. A derivative of order three is
not $L^2$, including any nonzero homogeneous combination of that order.
\end{lemma}
\begin{proof}
For $\phi\in C_c^\infty(\R^2\setminus\{0\})$ the sites eventually
avoid its support, so $R_h(\lambda)(\Delta^2+\lambda)\phi=\phi$.
The limiting operator therefore extends the free restriction. The range of its
resolvent difference lies in $\ker(S^*+\lambda)$ for that restriction.
Theorem~\ref{fd:thm:defect} with $m=n=2$ gives exactly the six
independent monomials of degree at most two and excludes degree three.
\end{proof}

\begin{lemma}[Limits of point resolvents]\label{lem:pseudo}
Let $H_h$ be self-adjoint, and suppose that for every nonreal $\lambda$
the resolvents $(H_h+\lambda)^{-1}$ converge in norm to
$R(\lambda)=R_0(\lambda)+\Psi_\lambda$, where $\ran\Psi_\lambda$
consists of functions $R_0(\lambda)T$ with $T$ a distribution supported
in a fixed finite set. Then $R(\lambda)=(H+\lambda)^{-1}$ for a
self-adjoint operator $H$.
\end{lemma}
\begin{proof}
The resolvent identity and $R(\lambda)^*=R(\bar\lambda)$ pass to the
limit, so $R$ is a pseudo-resolvent. If $R(\lambda)f=0$, then
$f=-(\Delta^2+\lambda)\Psi_\lambda f$ is an $L^2$ function whose
distributional support is finite, hence $f=0$. A pseudo-resolvent with trivial kernel
is the resolvent of a closed operator \cite[Chapter~VIII, Section~1, pp.~427--428]{Kato};
$\ker R(\bar\lambda)=0$ makes its domain dense, and
$R(\lambda)^*=R(\bar\lambda)$ makes it self-adjoint.
\end{proof}

\subsection{The critical collision limit}\label{sec:loc-critical}

Introduce the exact, energy-independent counterterm
\begin{equation}\label{eq:d}
 d_h=\frac{2(g_1-G_1(he))}{h^2}
 =\frac{\log(1/h)+1+\log2-\gamma_{\rm Eul}}{4\pi}
       +\frac{h^2}{256}+O(h^4(1+|\log h|)).
\end{equation}
The reference value $1$ fixes the unit and the coupling convention; changing
it shifts the scalar parameter below.

\begin{theorem}[Sufficient critical scaling]\label{thm:critical}
Let $a,\theta\in\R$. Suppose
\begin{equation}\label{eq:critical}
 \frac{b_1(h)}h\longrightarrow a,
 \qquad \frac{b_1(h)+b_2(h)}{h^2}+d_h\longrightarrow\theta.
\end{equation}
Then $H_h$ converges in norm resolvent sense to a self-adjoint operator
$H_{a,\theta}$. With $V_\lambda=(G_\lambda,D_\lambda)$, its resolvent is
\begin{align}
 (H_{a,\theta}+\lambda)^{-1}
 &=R_0(\lambda)-V_\lambda K_\lambda^{-1}V_{\bar\lambda}^*,\label{eq:limitres}\\
 K_\lambda&=\begin{pmatrix}
  \dfrac1{8\sqrt\lambda}&-a\\[3pt]
  -a&\theta-\dfrac{\log\lambda}{16\pi}
 \end{pmatrix}.\label{eq:K}
\end{align}
The statement holds at every nonreal $\lambda$; on the positive real axis
exclude the zeros of $\det K_\lambda$. The correction has rank two.
\end{theorem}
\begin{proof}
Write
\[
 C_h=\begin{pmatrix}1&1\\0&h\end{pmatrix},\quad
 C_h^{-1}=\begin{pmatrix}1&-h^{-1}\\0&h^{-1}\end{pmatrix},\quad
 V_{h,\lambda}=(G_\lambda,D_{h,\lambda}),\qquad
 \Phi_h=V_{h,\lambda}C_h.
\]
The transformed matrix $K_{h,\lambda}=C_h^{-T}Q_h(\lambda)C_h^{-1}$ is
\begin{equation}\label{eq:Kh}
 \begin{pmatrix}
  b_1+g_\lambda&[G_\lambda(he)-g_\lambda-b_1]/h\\
  [G_\lambda(he)-g_\lambda-b_1]/h&
  [b_1+b_2+2(g_\lambda-G_\lambda(he))]/h^2
 \end{pmatrix}.
\end{equation}
Equations \eqref{eq:critical} and Lemma~\ref{lem:kernel} show that
$K_{h,\lambda}\to K_\lambda$ and $V_{h,\lambda}\to V_\lambda$ in norm.
For $\operatorname{Im}\lambda>0$, both diagonal entries of
$\operatorname{Im}K_\lambda$ are strictly negative and the off-diagonal entry
is real. Thus $K_\lambda$ is invertible; the lower half-plane is analogous.
Substitution in \eqref{eq:res} proves norm convergence to \eqref{eq:limitres}.

By Lemma~\ref{lem:pseudo} the limit is the resolvent of a self-adjoint
operator. The two columns of $V_\lambda$ are independent; hence the
correction has rank two.
\end{proof}

\begin{remark}[A local limit is not just a distributional limit]
The sum $\alpha_1(h)\delta_0+\alpha_2(h)\delta_{he}$ need not converge as a
distribution. The object that converges is the associated operator resolvent.
The logarithmic counterterm in the \emph{inverse} strengths is essential.
\end{remark}

\subsection{The two-site converse}\label{sec:loc-converse}

This section concerns the fixed direction $e$ and the form class
\eqref{eq:form}; it is not a restriction on arbitrary regularizations.

\begin{theorem}[Necessary scaling and closure]\label{thm:converse}
Suppose $H_h$ has a norm resolvent limit which is a self-adjoint operator on
$L^2(\R^2)$. If its resolvent correction relative to $H_0$ has rank two, then
\eqref{eq:critical} holds for some real $a,\theta$ and the limit is
$H_{a,\theta}$. Every rank-one limit is a monopole point interaction at $0$;
the rank-zero limit is $H_0$. In particular, a dipole cannot occur alone in a
rank-one limit of these two scalar sites.
\end{theorem}
\begin{proof}
We first justify coefficient recovery under norm convergence of source
maps. If $W_h:\C^r\to L^2$ converges in norm to an injective $W$, its
least singular value is eventually at least half that of $W$. Hence
\begin{equation}\label{eq:left}
 \mathcal L_h=(W_h^*W_h)^{-1}W_h^*
       \longrightarrow(W^*W)^{-1}W^*=\mathcal L,
 \qquad \mathcal L_hW_h=I.
\end{equation}
This follows also directly from convergence of the positive Gram matrices
and their inverses. At a fixed nonreal $\lambda$, write
$B_h=K_{h,\lambda}^{-1}$. The limiting correction $R(\lambda)-R_0(\lambda)$
therefore determines
\begin{equation}\label{eq:recoverB}
 B_h=-\mathcal L_{h,\lambda}(R_h(\lambda)-R_0(\lambda))
                    \mathcal L_{h,\bar\lambda}^*\longrightarrow B.
\end{equation}
It follows that $R-R_0=-V_\lambda BV_{\bar\lambda}^*$ and that its rank
is the rank of $B$. Here $B^T=B$; Hermitian symmetry is not assumed at a
nonreal parameter.

\emph{Rank two.} If $B$ is invertible, $K_{h,\lambda}\to B^{-1}$.
The exact identity
\begin{equation}\label{eq:identity}
 (K_{h,\lambda})_{11}+h(K_{h,\lambda})_{12}=G_\lambda(he)
\end{equation}
forces $(B^{-1})_{11}=g_\lambda$. Since
$(G_\lambda(he)-g_\lambda)/h\to0$, its off-diagonal entry forces
$b_1/h\to a\in\R$. Its lower-right entry and \eqref{eq:renorm} force
$(b_1+b_2)/h^2+d_h\to\theta\in\R$. Theorem~\ref{thm:critical} then
identifies the limit.

\emph{All other matrix cases.} Suppose $B_{22}\ne0$ and put
$r_h=(B_h)_{12}/(B_h)_{22}\to r$. Inverting the $2\times2$ matrix gives
\[
 r_h=\frac{b_1+g_\lambda-G_\lambda(he)}{h(b_1+g_\lambda)},\qquad
 b_1+g_\lambda=\frac{G_\lambda(he)}{1-hr_h}.
\]
As $G_\lambda(he)-g_\lambda=O(h^2|\log h|)$, this implies
$b_1/h\to g_\lambda r\in\R$. The first row of $K_{h,\lambda}$
converges. The identity
$(K_h)_{12}(B_h)_{12}+(K_h)_{22}(B_h)_{22}=1$
then forces convergence of $(K_h)_{22}$. Thus $K_h\to K$ and $KB=I$;
this case has rank two. If $B_{22}=0$ and $B_{12}\ne0$, then
$\det B=-B_{12}^2\ne0$ but $(B^{-1})_{11}=0$, contradicting
\eqref{eq:identity}. Finally $B_{22}=B_{12}=0$ gives
$B=\diag(B_{11},0)$. It gives the free operator if $B_{11}=0$, and a
pure monopole correction otherwise. These cases exhaust the possibilities.

\emph{The scalar coupling is one real constant.} By
Lemma~\ref{lem:sixjets} a rank-one limit is a self-adjoint extension of
the six-jet restriction, and by the previous step its resolvent
correction has range $\operatorname{span}\{G_\lambda\}$. In
Theorem~\ref{fd:thm:trace}, with the jets and anchor of
Theorem~\ref{fd:thm:planarM}, this means that $\ran\Pi$ is the monopole
coordinate and $\Theta$ is a real number. Since
$(M_{-\lambda})_{11}=g_\lambda-g_1$ by \eqref{fd:eq:fullM}, formula
\eqref{fd:eq:projectionresolvent} reads
\begin{equation}\label{eq:monopole}
 R(\lambda)=R_0(\lambda)
       -G_\lambda(\beta+g_\lambda)^{-1}G_{\bar\lambda}^*,
 \qquad\beta=-\Theta-g_1\in\R.
\end{equation}
For $\beta\ne0$, this limit is attained with $b_1=\beta$ and
$b_2=h^{-2}$. For $\beta=0$, take $b_1=h^2$, $b_2=h^{-2}$.
Indeed direct inversion of \eqref{eq:Q} gives the coefficient limit
$\diag((\beta+g_\lambda)^{-1},0)$ in the site coordinates.
Taking both $b_j=h^{-2}$ gives the free limit. Exactly zero strengths
are handled by deleting channels; they cannot yield a rank-two correction
along a convergent subsequence.
\end{proof}

In the specific normalization \eqref{eq:K}, write
\begin{equation}\label{eq:Theta}
 K_\lambda=\Theta+\diag\left(g_\lambda,-\frac{\log\lambda}{16\pi}\right),
 \qquad \Theta=\begin{pmatrix}0&-a\\-a&\theta\end{pmatrix}.
\end{equation}
Thus two diagonal scalar couplings do not produce an arbitrary real symmetric
$2\times2$ matrix: $\Theta_{11}=0$ is forced. This statement is
coordinate dependent as written; the coordinate convention is fixed in
\eqref{eq:green} and \eqref{eq:d}. The family itself, and the absence of a
rank-one dipole limit, are operator statements. Hermitian coupling matrices with nonreal entries are excluded in this
fixed real basis: complex conjugation commutes with every approximant
and the relation $\overline{R_h(\lambda)\bar f}=R_h(\bar\lambda)f$
passes to the limit.

\subsection{A family with a common lower bound}\label{sec:loc-stable}

We now give a fixed prescription and prove its spectral lower bound and
its optimal relative energy coercivity.

\begin{theorem}[Stable prescription and coercivity]\label{thm:stable}
For $a,\theta\in\R$, set
\begin{equation}\label{eq:prescription}
 b_1(h)=ah+h^2,\qquad b_2(h)=-ah-h^2-h^2d_h+\theta h^2,
 \qquad\alpha_j(h)=b_j(h)^{-1}.
\end{equation}
For small $h$ the strengths are nonzero and have opposite signs.
For all sufficiently small $h$ the operators $H_h$ have a common
lower bound, and they converge to $H_{a,\theta}$. At every fixed sufficiently large $\lambda>0$, put
\[
 \nu_h(\lambda)=\inf_{u\in H^2\setminus\{0\}}
       \frac{q_h[u]+\lambda\norm{u}_2^2}{a_{0,\lambda}[u]}.
\]
Then
\begin{equation}\label{eq:optimalpair}
 \nu_h(\lambda)\sim
       \frac{-4\pi\det K_\lambda}{g_\lambda\log(1/h)}.
\end{equation}
In particular, for one fixed larger shift and constants $c,c'>0$,
\begin{equation}\label{eq:coercivity}
 q_h[u]+\lambda\norm{u}_2^2\ge c'\norm{u}_2^2+
                  \frac{c}{1+|\log h|}\norm{u}_{H^2}^2
\end{equation}
for every fixed $a$, including $a=0$.
\end{theorem}
\begin{proof}
For $a\ne0$, leading terms $\pm ah$ give the signs. For $a=0$,
$b_1=h^2>0$ and $b_2=h^2(\theta-1-d_h)<0$ eventually.
Choose $\lambda_0$ so large that $\theta-\log\lambda_0/(16\pi)<0$;
then $\det K_{\lambda_0}<0$ and
$\det Q_h(\lambda_0)=h^2\det K_{h,\lambda_0}<0$ for small $h$.

\emph{Schur complement and lower bound.} Label the positive and negative
sites $p,n$. Give $H^2$ the positive energy
$a_+[u]=a_{0,\lambda_0}[u]+b_p^{-1}|u(x_p)|^2$.
Subtracting the rank-one positive-site correction from the free evaluation
representer at $x_n$ gives squared norm
\[
 r_h=g_{\lambda_0}-\frac{G_{\lambda_0}(he)^2}{b_p+g_{\lambda_0}},
 \qquad |u(x_n)|^2\le r_h a_+[u].
\]
Its Schur complement satisfies the exact formula and limit
\begin{equation}\label{eq:pairS}
 S_h=b_n+r_h=\frac{\det Q_h(\lambda_0)}{b_p+g_{\lambda_0}},\qquad
 \frac{S_h}{h^2}\longrightarrow
              \frac{\det K_{\lambda_0}}{g_{\lambda_0}}<0.
\end{equation}
Hence
\[
 q_h[u]+\lambda_0\norm{u}_2^2
 =a_+[u]-|b_n|^{-1}|u(x_n)|^2
 \ge\frac{-S_h}{|b_n|}a_+[u]\ge0.
\]
This proves $H_h\ge-\lambda_0$. Since $a_+\ge a_{0,\lambda_0}$ and
$|b_n|\asymp h$ for $a\ne0$, $|b_n|\asymp h^2(1+|\log h|)$ for $a=0$,
it also gives the useful, weaker estimates
\begin{equation}\label{eq:coarseshur}
 q_h[u]+(\lambda_0+1)\norm{u}_2^2\ge\norm{u}_2^2+
 \begin{cases}ch\norm{u}_{H^2}^2,&a\ne0,\\
 c(1+|\log h|)^{-1}\norm{u}_{H^2}^2,&a=0.
 \end{cases}
\end{equation}
The scaling in Theorem~\ref{thm:critical} holds, since $b_1/h=a+h$.

\emph{Optimal relative constant.} Fix a sufficiently large $\lambda$ and
let $\tau:H^2\to\C^2$ be evaluation, with the energy $a_{0,\lambda}$.
Its Gram matrix is $M=M_\lambda>0$. The map $\tau^*M^{-1/2}$ is an
isometry onto the span of the two evaluation representers. In these
coordinates the relative form operator is $F_h=I+M^{1/2}\diag(b_1^{-1},b_2^{-1})M^{1/2};$ on the energy-orthogonal complement it is $I$. Its trace $t_h$ and
determinant $\delta_h$ are exactly
\[
 t_h=2+g_\lambda(b_1^{-1}+b_2^{-1}),\qquad
 \delta_h=\frac{\det Q_h}{b_1b_2}
             =\frac{h^2\det K_{h,\lambda}}{b_1b_2}.
\]
Writing $L=\log(1/h)$, for $a\ne0$ these have asymptotics
\[
 t_h\sim\frac{g_\lambda L}{4\pi a^2},\qquad
 \delta_h\longrightarrow-\frac{\det K_\lambda}{a^2}>0.
\]
For $a=0$, $t_h\sim\frac{g_\lambda}{h^2},\qquad
 \delta_h\sim-\frac{\det K_\lambda}{h^2d_h}.$

In both cases $\delta_h/t_h^2\to0$ and the smaller eigenvalue is
$\delta_h/t_h\,(1+o(1))$. This gives \eqref{eq:optimalpair}; it is
below $1$ and hence equals $\nu_h(\lambda)$. Increasing this fixed shift
by $1$ gives \eqref{eq:coercivity} with a separate $L^2$ term.
\end{proof}

\begin{remark}[Slow onset of the asymptotic regime]\label{rem:slowonset}
For $a\ne0$ the regime described by \eqref{eq:optimalpair} is reached
slowly. With $t_h$ and $\delta_h$ as in the preceding proof, the exact
value for small $h$ is
$\nu_h(\lambda)=\tfrac12\bigl(t_h-\sqrt{t_h^2-4\delta_h}\bigr)$, where
\[
 t_h=2+\frac{g_\lambda(d_h-\theta)}{a^2}\bigl(1+O(h\log(1/h))\bigr),
\]
and the logarithmic term dominates the constant $2$ only when
$\log(1/h)\gg64\pi a^2\sqrt\lambda$.
\end{remark}

\section{Smooth local potentials and quantitative pair limits}\label{sec:loc-smoothing}

\subsection{Realization by smooth local potentials}\label{sec:loc-smooth}

Fix $\rho\in C_c^\infty(\R^2)$, $\supp\rho\subset\overline{B(0,1)}$,
$\rho\ge0$, $\int\rho=1$, and put
$\rho_\varepsilon(x)=\varepsilon^{-2}\rho(x/\varepsilon)$. This single
profile is fixed throughout; neither radiality nor a vanishing first
moment is assumed. Set $L_h=1+|\log h|$, $L_\varepsilon=1+|\log\varepsilon|$.

\begin{theorem}[Smooth local realization]\label{thm:smooth}
Use the strengths \eqref{eq:prescription} and define
\begin{equation}\label{eq:potential}
 W_{h,\varepsilon}(x)=\alpha_1(h)\rho_\varepsilon(x)
                 +\alpha_2(h)\rho_\varepsilon(x-he),\qquad
 \widetilde H_{h,\varepsilon}=\Delta^2+W_{h,\varepsilon},\quad
 \Dom(\widetilde H_{h,\varepsilon})=H^4(\R^2).
\end{equation}
For $\varepsilon_h=h^p$ with any fixed $p>1$, these operators have a
common lower bound and converge in norm resolvent sense to $H_{a,\theta}$.
At a fixed sufficiently large real $\lambda$, for small $h$,
\begin{equation}\label{eq:smootherr}
 \norm{(\widetilde H_{h,\varepsilon_h}+\lambda)^{-1}-(H_h+\lambda)^{-1}}
 \le C\left[\varepsilon_h+
                \frac{\varepsilon_h^2L_{\varepsilon_h}L_h}{h^2}\right].
\end{equation}
The constants may depend on $p,a,\theta,\lambda,\rho$, but not on $h$.
In particular $\varepsilon_h=h^2$ suffices. The choice $h^8$ gives error $O(h^8)$.
The limiting operator domain is not contained in $H^2(\R^2)$.
\end{theorem}
\begin{proof}
\emph{Modulus of continuity and trace estimates.}
Fourier Cauchy--Schwarz and a split at $|\xi|=|y|^{-1}$ give, for
$0<|y|<1/2$,
\begin{align}
 |u(x+y)-u(x)|^2
 &\le C\norm{u}_{H^2}^2\int_{\R^2}
              \frac{|e^{i\xi\cdot y}-1|^2}{1+|\xi|^4}\,\dd\xi
 \le C|y|^2(1+|\log|y||)\norm{u}_{H^2}^2.\label{eq:modulus}
\end{align}
The low-frequency part is bounded by
$C|y|^2\int_0^{1/|y|}r^3(1+r^4)^{-1}\dd r$ and the high-frequency
part by $C\int_{1/|y|}^\infty r(1+r^4)^{-1}\dd r$.
Also $\norm{u}_\infty\le C\norm{u}_{H^2}$. These imply the global
$C^{0,1/2}$ estimate and, uniformly in the center $x_j$,
\begin{equation}\label{eq:mollify}
 \left|\int\rho_\varepsilon(x-x_j)|u(x)|^2\dd x-|u(x_j)|^2\right|
 \le C\varepsilon^{1/2}\norm{u}_{H^2}^2.
\end{equation}
Here we used $|y|\sqrt{1+|\log|y||}\le C|y|^{1/2}$ for small $|y|$.

\emph{Averaged traces and their Gram matrix.} For the moment take any
finite site set with the same profile and define
\[
 \ell_{\varepsilon,j}(u)=\int\rho_\varepsilon(x-x_j)u(x)\dd x,\qquad
 \overline q_{h,\varepsilon}[u]=\norm{\Delta u}_2^2+
                          \sum_j\alpha_j|\ell_{\varepsilon,j}(u)|^2.
\]
This intermediate form is used only in the proof. Its source columns are
$P_\varepsilon G_\lambda(\cdot-x_j)$, where
$P_\varepsilon f=\rho_\varepsilon*f$, and its free Gram kernel is
\[
 M_{\lambda,\varepsilon}(x)=(2\pi)^{-2}\int_{\R^2}
 \frac{e^{ix\cdot\xi}|\widehat\rho(\varepsilon\xi)|^2}
      {|\xi|^4+\lambda}\,\dd\xi.
\]
Nonnegativity and unit mass of $\rho$ give
\[
 0\le1-|\widehat\rho(\varepsilon\xi)|^2
 =\iint\rho(s)\rho(t)[1-\cos(\varepsilon\xi\cdot(s-t))]\dd s\dd t
 \le C\min(\varepsilon^2|\xi|^2,1).
\]
At a fixed positive shift this proves the uniform estimate
\begin{equation}\label{eq:meanGram}
 \sup_x|M_{\lambda,\varepsilon}(x)-G_\lambda(x)|
 \le C\left[\varepsilon^2\int_0^{1/\varepsilon}
                   \frac{r^3}{1+r^4}\dd r
            +\int_{1/\varepsilon}^\infty\frac{r}{1+r^4}\dd r\right]
 \le C\varepsilon^2L_\varepsilon.
\end{equation}
The kernel is even since $|\widehat\rho|^2$ is even. Independence
of the averaged traces follows because $\widehat\rho$ is nonzero near
zero and a finite sum of distinct site exponentials vanishing there
vanishes identically. Hence their Gram matrix is positive definite at
a positive shift, and \eqref{eq:res} applies to the averaged traces.

For a pair, $\norm{C_h^{-1}}=O(h^{-1})$, so the transformed matrix changes
by $O(\varepsilon^2L_\varepsilon/h^2)=o(1)$ for $p>1$.
Its source map is $P_\varepsilon V_{h,\lambda}$. Since $P_\varepsilon$
is an $L^2$ contraction, tends strongly to $I$, and the source space is
finite dimensional, $P_\varepsilon V_{h,\lambda}\to V_\lambda$ in norm.
Thus the averaged resolvents have the same limit.

Their lower bound follows directly from the same Schur calculation.
Writing $m_0=M_{\lambda_0,\varepsilon}(0)$ and
$m_h=M_{\lambda_0,\varepsilon}(he)$, the positive-site denominator is
$b_p+m_0\to g_{\lambda_0}$ and
\[
 S_{h,\varepsilon}
 =b_n+m_0-\frac{m_h^2}{b_p+m_0}
 =\frac{h^2\det\overline K_{h,\lambda_0}}{b_p+m_0},\qquad
 S_{h,\varepsilon}/h^2\to\det K_{\lambda_0}/g_{\lambda_0}<0.
\]
Here $\overline K$ denotes the transformed averaged matrix. The positive
averaged-site energy dominates $a_{0,\lambda_0}$, so after a fixed shift
increase the averaged forms obey \eqref{eq:coarseshur} with uniform
constants, and in particular a lower bound by $c'I$.

\emph{Return to the local multiplication form.} The exact variance
identity is
\begin{align}
 &\int\rho_\varepsilon(x-x_j)|u(x)|^2\dd x
                 -|\ell_{\varepsilon,j}(u)|^2\nonumber\\
 &\quad=\tfrac12\iint\rho(s)\rho(t)
          |u(x_j+\varepsilon s)-u(x_j+\varepsilon t)|^2\dd s\dd t
 \le C\varepsilon^2L_\varepsilon\norm{u}_{H^2}^2.\label{eq:variance}
\end{align}
The last bound follows from \eqref{eq:modulus}; the factor
$|s-t|^2|\log|s-t||$ is bounded on the fixed compact integration set.
For $a\ne0$, $\sum|\alpha_j|=O(h^{-1})$; for $a=0$ it is $O(h^{-2})$.
Using the averaged Schur estimates just proved, the signed sum of the
variances has relative shifted-form bound
\begin{equation}\label{eq:varianceeta}
 \eta_{h,\varepsilon}\le C
 \begin{cases}
 \varepsilon^2L_\varepsilon/h^2,&a\ne0,\\
 \varepsilon^2L_\varepsilon L_h/h^2,&a=0.
 \end{cases}
\end{equation}
These bounds tend to zero for $p>1$.

We use the following form perturbation estimate. Let $A\ge c'I$ have
closed form $a$, and let $t$ be a Hermitian form on its form domain with
$|t[u]|\le\eta a[u]$, $\eta<1$. Then $a+t$ is closed and positive
\cite[Section~6.5]{Teschl}; writing $t(u,v)=\ip{A^{1/2}u}{TA^{1/2}v}$ with
$T=T^*$, $\norm T\le\eta$, its operator has inverse
$A^{-1/2}(I+T)^{-1}A^{-1/2}$. Hence
\begin{equation}\label{eq:formpert}
 \norm{R_{a+t}-A^{-1}}\le(c')^{-1}\frac{\eta}{1-\eta},
\end{equation}
and the shifted local form is positive, which gives a common lower bound.

\emph{Quantitative comparison and the limiting domain.}
For $f\in H^1$, translation and integration against $\rho$ give
$\norm{(P_\varepsilon-I)f}_2\le C\varepsilon\norm f_{H^1}$.
The pair source columns are uniformly bounded in $H^1$, since
$|e^{-ih\xi_e}-1|/h\le|\xi_e|$. Hence their change under averaging is
$O(\varepsilon)$. The transformed point and averaged matrix inverses
are uniformly bounded at the chosen shift. Subtracting their resolvent
formulas using \eqref{eq:meanGram} gives
$O(\varepsilon+\varepsilon^2L_\varepsilon/h^2)$.
Add \eqref{eq:formpert} with \eqref{eq:varianceeta} to get
\eqref{eq:smootherr}. For fixed $h$, the local potential is a bounded
smooth real multiplication operator, so its operator domain is $H^4$.

Finally $V_{\bar\lambda}^*$ is onto and $K_\lambda$ is invertible.
For a suitable $f\in L^2$, the resolvent formula has a nonzero
$D_\lambda$ coefficient. Its other terms are in $H^2$, whereas
$D_\lambda\notin H^2$, proving the domain assertion.
\end{proof}

\subsubsection{The filter loss and the critical smoothing scale}
Let $m_\rho=\int x\rho(x)\dd x$ and put
\[
 \Sigma_\rho=\int(x-m_\rho)(x-m_\rho)^T\rho(x)\dd x,
 \qquad \kappa_\rho=\frac{\tr\Sigma_\rho}{2\pi}>0.
\]
This covariance is positive definite: a zero variance in a direction
would put the support of $\rho$ in a line, whereas a smooth nonnegative
unit-mass function is positive on a nonempty open set. No centering or
radiality assumption is added.

\begin{lemma}[Loss of a divided-difference Gram entry]\label{lem:filter}
For fixed $\lambda>0$ set
\begin{align*}
 I_1(h,\varepsilon)&=(2\pi)^{-2}\int
 \frac{(1-|\widehat\rho(\varepsilon\xi)|^2)
                         |e^{-ih\xi_e}-1|^2}
 {h^2(|\xi|^4+\lambda)}\dd\xi,\\
 I_2(h,\varepsilon)&=(2\pi)^{-2}\int
 \frac{(1-|\widehat\rho(\varepsilon\xi)|^2)
                    |e^{-ih\xi_e}+e^{ih\xi_e}-2|^2}
 {h^4(|\xi|^4+\lambda)}\dd\xi.
\end{align*}
If $t=\varepsilon/h\to0$, then, uniformly for small $h$,
\begin{align}
 I_1(h,\varepsilon)&=\kappa_\rho t^2\log(1/t)+O(t^2),\label{eq:filter1}\\
 I_2(h,\varepsilon)&=3\kappa_\rho\frac{\varepsilon^2}{h^4}
                         \log(h/\varepsilon)
                           +O(\varepsilon^2/h^4).\label{eq:filter2}
\end{align}
For any fixed $r>0$, $I_1(h,rh)$ tends to the strictly positive constant
\begin{equation}\label{eq:Cprofile}
 C_\rho(r,e)=(2\pi)^{-2}\int
 \frac{(1-|\widehat\rho(r\eta)|^2)|e^{-i\eta_e}-1|^2}{|\eta|^4}
 \dd\eta.
\end{equation}
If $\varepsilon=h^p$, $0<p<1$, then
\begin{equation}\label{eq:widepairGram}
 \frac{2[M_{\lambda,\varepsilon}(0)-M_{\lambda,\varepsilon}(he)]}{h^2}
                   =\frac p{4\pi}\log(1/h)+O(1).
\end{equation}
\end{lemma}
\begin{proof}
Taylor expansion of the cosine in the double integral defining
$|\widehat\rho|^2$ gives
\[
 1-|\widehat\rho(\zeta)|^2=\zeta^T\Sigma_\rho\zeta+O(|\zeta|^4).
\]
The remainder is fourth order because the integrand is even in $\zeta$.
Also this difference is nonnegative and bounded by
$C\min(|\zeta|^2,1)$. Put
$p_1(s)=2-2\cos s$ and $p_2(s)=6-8\cos s+2\cos(2s)$.
The substitution $\eta=h\xi$ gives the prefactor $1$ for $I_1$
and $h^{-2}$ for $I_2$, with denominator $|\eta|^4+\lambda h^4$.
On $|\eta|<1$, use $p_j(\eta_e)=O(|\eta|^{2j})$ and the quadratic
bound on the filter. This part is $O(t^2)$ before the prefactor.
On $|\eta|>1/t$ its absolute value is also $O(t^2)$.
On $1<|\eta|<1/t$, replacing the filter by
$t^2\eta^T\Sigma_\rho\eta$ costs $O(t^2)$, since
$t^4\int_1^{1/t}s\dd s=O(t^2)$.
Replacing the denominator there by $|\eta|^4$ has the same bound.

The constant terms of $p_1,p_2$ are $2,6$, respectively, and
$\int_{\mathbb S^1}\omega^T\Sigma_\rho\omega\dd\omega
=\pi\tr\Sigma_\rho$. They give the coefficients in
\eqref{eq:filter1}--\eqref{eq:filter2}.
Their oscillatory terms contribute only $O(t^2)$: for $k=1,2$,
\[
 \left|\int_1^R\frac{\cos(ks\omega_e)}s\dd s\right|
        \le C+|\log|\omega_e||.
\]
Split the integral at $1/|\omega_e|$ and integrate by parts beyond
that point to obtain this inequality. Its right side is integrable
on the circle, so angular integration is uniformly bounded.

For \eqref{eq:Cprofile}, use the same substitution with $t=r$ and
dominated convergence. Near zero its numerator is $O(|\eta|^4)$;
at infinity the denominator gives an integrable bound. The same bounds
are uniform for $t$ in a compact subset of $(0,\infty)$, so the conclusion
also holds when $t\to r$. The integral
is positive: equality $|\widehat\rho(\zeta)|=1$ for a nonzero
$\zeta$ would force $e^{-i\zeta\cdot x}$ to have constant phase on
an open set where $\rho>0$, which is impossible.

Finally, if $h/\varepsilon\to0$, Taylor expansion of the first stencil
gives
\[
 \frac{2(M_{\lambda,\varepsilon}(0)-M_{\lambda,\varepsilon}(he))}{h^2}
 =(2\pi)^{-2}\int\frac{\xi_e^2|\widehat\rho(\varepsilon\xi)|^2}
 {|\xi|^4+\lambda}\dd\xi+O(h^2/\varepsilon^2).
\]
The remainder follows by bounding the fourth moment after the Fourier
cutoff. To see its finiteness, integration by parts in the transform
gives $|\widehat\rho(\zeta)|\le C_N(1+|\zeta|)^{-N}$.
Split the main integral at $1$ and $1/\varepsilon$.
On the middle interval replacing $|\widehat\rho|^2$ by $1$ costs
$O(\varepsilon^2\int_1^{1/\varepsilon}s\dd s)=O(1)$;
the tail is $O(1)$ by this rapid decay. The angular integral of
$\omega_e^2$ is $\pi$, giving
$\log(1/\varepsilon)/(4\pi)+O(1)$.
Substitution $\varepsilon=h^p$ proves \eqref{eq:widepairGram}.
\end{proof}

\begin{theorem}[Sharp pair powers and the comparable-width limit]
\label{thm:pairwidth}
Fix $a\ne0$ and $\theta$, and keep exactly the inverse strengths
\eqref{eq:prescription}. For $\varepsilon=h^p$ with any $p>0$, the
smooth local operators have a common lower bound. Their limits are
\begin{equation}\label{eq:pairphases}
 \begin{cases}
 H_{a,\theta},&p>1,\\
 H_{a,\theta-C_\rho(1,e)},&p=1,\\
 H_{\rm pin},&0<p<1,
 \end{cases}
\end{equation}
where
\begin{equation}\label{eq:Rpin}
 (H_{\rm pin}+\lambda)^{-1}
   =R_0(\lambda)-G_\lambda g_\lambda^{-1}G_{\bar\lambda}^*.
\end{equation}
The form domain of $H_{\rm pin}$ is
$\{u\in H^2(\R^2):u(0)=0\}$.
More generally $\varepsilon/h\to r\in(0,\infty)$ gives
$H_{a,\theta-C_\rho(r,e)}$.
For $a\ne0$, arbitrary separated widths $\varepsilon/h\to0$ give
$H_{a,\theta}$. Consequently $p>1$ is the exact power range preserving
the originally prescribed limit for fixed nonzero $a$ and fixed profile.
The comparable and wider regimes at $a=0$ are not asserted here.
\end{theorem}
\begin{proof}
Denote averaged quantities by a bar. Exactly,
\[
 \overline K_{11}=b_1+m_0,\quad
 \overline K_{12}=-b_1/h-\frac h2\overline T_h,\quad
 \overline K_{22}=(b_1+b_2)/h^2+\overline T_h,
\]
where $m_0=M_{\lambda,\varepsilon}(0)$ and
$\overline T_h=2(M_{\lambda,\varepsilon}(0)
-M_{\lambda,\varepsilon}(he))/h^2$.
The difference between its last entry and the point entry is exactly
$-I_1(h,\varepsilon)$; its off-diagonal difference is $hI_1/2$.
For a finite positive width ratio, Lemma~\ref{lem:filter} therefore
gives $\overline K\to K_\lambda(a,\theta-C_\rho(r,e))$.
For a ratio tending to zero it gives $\overline K\to K_\lambda(a,\theta)$.
The source maps $P_\varepsilon V_{h,\lambda}$ converge to $V_\lambda$,
so these are the averaged norm resolvent limits.

For $0<p<1$, \eqref{eq:widepairGram} instead gives
\[
 \overline K_{11}\to g_\lambda,\quad
 \overline K_{12}\to-a,\quad
 \overline K_{22}=-\frac{1-p}{4\pi}\log(1/h)+O(1).
\]
Thus $\overline K^{-1}\to\diag(g_\lambda^{-1},0)$ and the averaged
limit is \eqref{eq:Rpin}. The kernel of evaluation is a closed subspace
of $H^2$ in the free shifted energy and is dense in $L^2$, since it
contains smooth functions supported away from zero.
For $u_0=R_0(\lambda)f$, the energy projection onto this subspace is
$u_0-G_\lambda u_0(0)/g_\lambda$. It has value zero at the origin;
pairing its free shifted energy with a test vector of value zero removes
the $G_\lambda$ term and gives the required weak equation. This proves
\eqref{eq:Rpin} and the asserted form domain.

We prove the needed bounds before passing to multiplication potentials.
At one fixed sufficiently large shift the averaged relative form
matrix has trace and determinant
\[
 t_h=2+m_0(b_1^{-1}+b_2^{-1}),\qquad
 \delta_h=h^2\det\overline K/(b_1b_2).
\]
For fixed $a\ne0$, $b_1b_2\sim-a^2h^2$ and
$b_1^{-1}+b_2^{-1}\sim\log(1/h)/(4\pi a^2)$.
In the finite-ratio or separated regimes, choose the shift so that
the limiting determinant is negative. Then the smaller relative
eigenvalue is asymptotic to $\delta_h/t_h$ and is bounded below by
$c/\log(1/h)$. For $0<p<1$ it tends to $1-p$:
\[
 t_h\sim\frac{g_\lambda L}{4\pi a^2},\qquad
 \delta_h\sim\frac{g_\lambda(1-p)L}{4\pi a^2},\quad L=\log(1/h).
\]
Here $\delta_h/t_h^2\to0$; the complement of the evaluation span
has relative value $1$. These estimates prove common averaged lower
bounds. Increase the fixed shift by one to retain a separate positive
$L^2$ term for the form perturbation argument.

For separated widths the unsigned variance estimate now uses this
optimal averaged coercivity and gives
\begin{equation}\label{eq:pairvariancebetter}
 \eta_{h,\varepsilon}\le
             C\varepsilon^2L_\varepsilon L_h/h=o(1).
\end{equation}
Indeed, writing $t=\varepsilon/h$, its right side is bounded by
$Ch[t^2L_h^2+t^2\log(1/t)L_h]$, which tends to zero.
For a finite positive width ratio the same bound is $O(hL_h^2)$.
Thus \eqref{eq:formpert} proves the claimed local limits and bounds
in these regimes.

For wider powers, use cancellation of the signed variances rather
than their absolute sum. If $\mathcal V_\varepsilon(x;u)$ denotes
the variance in \eqref{eq:variance}, the modulus estimate gives
\[
 |\mathcal V_\varepsilon(x+he;u)-\mathcal V_\varepsilon(x;u)|
       \le Ch\varepsilon\sqrt{L_hL_\varepsilon}\norm u_{H^2}^2.
\]
To verify it, subtract the squared increments in the double integral.
Their sum has bound $C\varepsilon\sqrt{L_\varepsilon}\norm u_{H^2}$,
and their difference has bound $Ch\sqrt{L_h}\norm u_{H^2}$.
Now $\alpha_1+\alpha_2=O(L_h)$ and $\alpha_2=O(h^{-1})$.
Their signed variance sum is bounded in absolute value by
\[
 C[L_h\varepsilon^2L_\varepsilon+
                 \varepsilon\sqrt{L_hL_\varepsilon}]\norm u_{H^2}^2.
\]
For every fixed $0<p<1$ the bracket tends to zero and the averaged
coercivity stays positive. Equation~\eqref{eq:formpert} applies again.
Since $C_\rho(1,e)>0$, the $p=1$ limit is different from the original
rank-two limit. The pinned limit has rank-one correction and is
different as well. This proves the sharp power assertion.
\end{proof}

\subsection{The leading resolvent error}\label{sec:loc-rate}

The following is a rate statement for \eqref{eq:prescription}, not an optimal
rate among every possible tuning of the strengths.

\begin{theorem}[Leading collision error]\label{thm:rate}
Let $\lambda>0$ be a fixed sufficiently large shift and put
$J=\begin{psmallmatrix}0&1\\1&0\end{psmallmatrix}$. Both the point model and
its smooth realization satisfy the operator norm limit
\begin{equation}\label{eq:leading}
 \frac{(H_h+\lambda)^{-1}-(H_{a,\theta}+\lambda)^{-1}}{h\log h}
 \longrightarrow \frac1{8\pi}\mathcal J_\lambda,
 \qquad
 \mathcal J_\lambda:=V_\lambda K_\lambda^{-1}J K_\lambda^{-1}V_\lambda^*.
\end{equation}
For the smooth model take any fixed $p>3/2$ in Theorem~\ref{thm:smooth}
and replace $H_h$ by $\widetilde H_{h,h^p}$ in the numerator.
The bounded self-adjoint operator $\mathcal J_\lambda$ is nonzero. In particular,
the resolvent error is asymptotic to
$\norm{\mathcal J_\lambda}h|\log h|/(8\pi)$.
\end{theorem}
\begin{proof}
From \eqref{eq:Kh}, \eqref{eq:prescription}, and Lemma~\ref{lem:kernel},
\[
 K_{h,\lambda}=K_\lambda+\frac{h\log h}{8\pi}J+O(h),
 \qquad V_{h,\lambda}=V_\lambda+O(h).
\]
Expand the inverse of $K_{h,\lambda}$ to first order and substitute in the
resolvent formula. The quadratic remainder is
$O(h^2|\log h|^2)=o(h)$, giving \eqref{eq:leading}. Let $\mathcal L=(V_\lambda^*V_\lambda)^{-1}V_\lambda^*$.
If $\mathcal J_\lambda$ vanished, multiplying by $\mathcal L$ and
$\mathcal L^*$ would give $K_\lambda^{-1}JK_\lambda^{-1}=0$.
This is impossible since both $K_\lambda$ and $J$ are invertible. For $\varepsilon=h^p$ with $p>3/2$, the right side of
\eqref{eq:smootherr} is
$O(h^p+h^{2p-2}L_h^2)=o(h|\log h|)$, so the smooth model has the same
leading operator. Norm continuity gives the last assertion.
\end{proof}

\begin{theorem}[The smoothing contribution to the leading pair error]
\label{thm:pairtwoscale}
Fix $a\ne0$ and use \eqref{eq:prescription}. Let
$\varepsilon/h\to0$, $L=\log(1/h)$, and
$\chi_h=(\varepsilon/h)^2\log(h/\varepsilon)$.
At a fixed sufficiently large real shift, let $\mathcal J_\lambda$ be
as in Theorem~\ref{thm:rate}, set $P_2=\diag(0,1)$ and
\[
 \mathcal P_\lambda=V_\lambda K_\lambda^{-1}P_2K_\lambda^{-1}V_\lambda^*.
\]
Then
\begin{equation}\label{eq:pairtwoscale}
 (\widetilde H_{h,\varepsilon}+\lambda)^{-1}
 -(H_{a,\theta}+\lambda)^{-1}
 =-\frac{hL}{8\pi}\mathcal J_\lambda
      -\kappa_\rho\chi_h\mathcal P_\lambda+o(hL+\chi_h).
\end{equation}
The norm error is comparable to $hL+\chi_h$. For $\varepsilon=h^p$,
$1<p<3/2$, division by $h^{2p-2}L$ gives the nonzero rank-one limit
$-\kappa_\rho(p-1)\mathcal P_\lambda$.
At $p=3/2$, division by $hL$ gives
\[
 -V_\lambda K_\lambda^{-1}
       \left[\frac{J}{8\pi}+\frac{\kappa_\rho}{2}P_2\right]
                                  K_\lambda^{-1}V_\lambda^*.
\]
Thus the strict threshold $p>3/2$ for retaining the leading operator
of Theorem~\ref{thm:rate} is sharp for fixed nonzero $a$.
\end{theorem}
\begin{proof}
The point matrix differs from $K_\lambda$ by
$-hLJ/(8\pi)+O(h)$. Lemma~\ref{lem:filter} changes its $(2,2)$
entry by $-\kappa_\rho\chi_h+o(\chi_h)$, its off-diagonal by
$O(h\chi_h)$ and its first diagonal by $O(\varepsilon^2L_\varepsilon)$.
The last quantity is $o(\chi_h)$: the ratio is bounded by
$Ch^2[1+L/\log(h/\varepsilon)]\to0$.
The source changes are $O(h+\varepsilon)=o(hL+\chi_h)$.
Finally \eqref{eq:pairvariancebetter} is $o(\chi_h)$, since its ratio
to $\chi_h$ is bounded by
$ChL[1+L/\log(h/\varepsilon)]\to0$.
Expand the bounded inverse matrix to first order. The quadratic
remainder is $o(hL+\chi_h)$, and \eqref{eq:formpert} transfers the
averaged expansion to the local multiplication operator. This proves
\eqref{eq:pairtwoscale} with the indicated signs.

The map $X\mapsto V_\lambda K_\lambda^{-1}XK_\lambda^{-1}V_\lambda^*$
is injective on $2\times2$ matrices and has a bounded inverse on its
finite-dimensional range. Matrices $J/(8\pi)+s\kappa_\rho P_2$,
$s\ge0$, have norms bounded below by a positive multiple of $1+s$:
their fixed off-diagonal is nonzero and their lower-right entry grows
linearly. This proves comparability and prevents cancellation.
For powers, $\chi_h=(p-1)h^{2p-2}L$; the asserted limits follow.
\end{proof}

\subsection{Spectral and sign consequences}\label{sec:loc-spectrum}

\begin{proposition}[Exactly one negative eigenvalue of the limit]\label{prop:pairnegative}
$H_{a,\theta}$ has essential spectrum $[0,\infty)$ and exactly one negative
eigenvalue $-\lambda_*$, which is simple. It is determined by
\begin{equation}\label{eq:eig}
 \theta-\frac{\log\lambda_*}{16\pi}=8a^2\sqrt{\lambda_*}.
\end{equation}
When $a=0$, $\lambda_*=e^{16\pi\theta}$.
\end{proposition}
\begin{proof}
Dividing $\det K_\lambda=0$ by $g_\lambda>0$ gives \eqref{eq:eig}.
The left side minus the right side decreases strictly from $+\infty$ to
$-\infty$ on $(0,\infty)$, with derivative
$-1/(16\pi\lambda)-4a^2/\sqrt\lambda$. Thus the determinant has one
simple zero. There $K_\lambda$ has rank one, and \eqref{eq:limitres}
has a nonzero rank-one residue since $V_\lambda$ is injective. This
gives the stated multiplicity. The finite-rank resolvent difference
preserves the essential spectrum \cite[Theorem~6.20]{Teschl}.
\end{proof}

\begin{proposition}[A sign restriction under a common lower bound]
A rank-two dipole limit from two scalar sites with a common lower bound
cannot be obtained from two repulsive strengths, or from two attractive
strengths. In any such convergent stable family, the two signs must be
opposite for all sufficiently small $h$.
\end{proposition}
\begin{proof}
Theorem~\ref{thm:converse} gives $b_1+b_2=h^2(\theta-d_h)+o(h^2)<0$.
At any sufficiently large fixed shift above the common lower bound in
absolute value, $\det Q_h=h^2\det K_{h,\lambda}<0$ eventually.
Theorem~\ref{thm:inertia} then gives $n_-(B_h)=n_-(Q_h)=1$.
A zero strength would leave at most one site and a correction of rank
at most one.
\end{proof}

\section{Three-site Hessian limits}\label{sec:loc-hessian}

Put the sites at $-he,0,he$, with equal inverse strengths $b_o$ at the
outer sites and $b_c$ at the center. For $h$ such that $b_o(h)b_c(h)\ne0$, let $H_h^{(3)}$ be associated with
\[
 q_h^{(3)}[u]=\norm{\Delta u}_2^2+
 b_o^{-1}(|u(-he)|^2+|u(he)|^2)+b_c^{-1}|u(0)|^2,
 \qquad\Dom(q_h^{(3)})=H^2(\R^2).
\]
The tuning below ensures these inverses exist for every sufficiently
small $h$. Recall $E_\lambda=\partial_e^2G_\lambda\in L^2\setminus H^1$.
For the fixed nonnegative unit-mass profile $\rho$ of
Section~\ref{sec:loc-smooth}, put
\[
 \widetilde H_{h,\varepsilon}^{(3)}=\Delta^2+
 b_o(h)^{-1}\bigl(\rho_\varepsilon(\cdot+he)
                    +\rho_\varepsilon(\cdot-he)\bigr)
                  +b_c(h)^{-1}\rho_\varepsilon,
 \qquad \Dom(\widetilde H_{h,\varepsilon}^{(3)})=H^4(\R^2).
\]
For each such $h$ and $\varepsilon>0$ this is a bounded real
multiplication perturbation of $\Delta^2$, with form domain $H^2$.

\begin{theorem}[Stable three-site Hessian limit]\label{thm:triplet}
For $\beta,\vartheta\in\R$, prescribe
\begin{align}
 b_c(h)&=G_1(he)-g_1+\beta h^2,\nonumber\\
 b_o(h)&=-2b_c(h)-3g_1+4G_1(he)-G_1(2he)
                            +\tfrac12\vartheta h^4.\label{eq:triscaling}
\end{align}
For small $h$, $b_o>0>b_c$. The operators $H_h^{(3)}$ have a common
lower bound and converge in norm resolvent sense to a self-adjoint
operator $H_{\beta,\vartheta}^{(E)}$, with resolvent
\begin{align}
 (H_{\beta,\vartheta}^{(E)}+\lambda)^{-1}
 &=R_0(\lambda)-U_\lambda A_\lambda^{-1}U_{\bar\lambda}^*,
 \qquad U_\lambda=(G_\lambda,E_\lambda),\label{eq:Eres}\\
 A_\lambda&=\begin{pmatrix}
 g_\lambda&-2\beta+\dfrac{\log\lambda}{16\pi}\\[3pt]
 -2\beta+\dfrac{\log\lambda}{16\pi}&
 \vartheta+\dfrac3{64}(1-\sqrt\lambda)
 \end{pmatrix}.\label{eq:EA}
\end{align}
This formula holds for every nonreal $\lambda$ and, on the positive real
axis, wherever $\det A_\lambda\ne0$.
At a fixed sufficiently large real shift $\lambda$,
\begin{equation}\label{eq:trirate}
 \log(1/h)\bigl[(H_h^{(3)}+\lambda)^{-1}
              -(H_{\beta,\vartheta}^{(E)}+\lambda)^{-1}\bigr]
 \longrightarrow-\frac{8\pi}{3}D_\lambda D_\lambda^*.
\end{equation}
The norm error is consequently asymptotic to
$[6\lambda\log(1/h)]^{-1}$. With $c_0=1/(8\pi)$, the optimal relative
energy constant has the asymptotic
\begin{equation}\label{eq:optimaltri}
 \inf_{u\in H^2\setminus\{0\}}
 \frac{q_h^{(3)}[u]+\lambda\norm u_2^2}{a_{0,\lambda}[u]}
 \sim\frac{-\det A_\lambda}{8c_0g_\lambda\log2}\,h^2.
\end{equation}
The limiting operator domain is not contained in $H^1(\R^2)$.
\end{theorem}
\begin{proof}
\emph{Source basis and the full congruence.} Define
\[
 D_{h,\lambda}^{s}=
 \frac{G_\lambda(\cdot-he)-G_\lambda(\cdot+he)}{2h},\qquad
 E_{h,\lambda}=
 \frac{G_\lambda(\cdot-he)+G_\lambda(\cdot+he)-2G_\lambda}{h^2}.
\]
The three source columns $(G_\lambda,D_{h,\lambda}^{s},E_{h,\lambda})$
converge in $L^2$ to $(G_\lambda,D_\lambda,E_\lambda)$ with respective
errors $0,O(h),O(h)$. The first derivative estimate follows from
translation Taylor. For the Hessian error its numerator, after subtracting
the limiting multiplier, is bounded by
$C\min(h^2|\xi|^4,|\xi|^2)$; hence its squared norm is at most
\[
 C\left[h^4\int_0^{1/h}\frac{r^9}{(1+r^4)^2}\dd r
       +\int_{1/h}^\infty\frac{r^5}{(1+r^4)^2}\dd r\right]=O(h^2).
\]
In the site order $(-he,0,he)$ the exact source change is
\begin{equation}\label{eq:C3}
 C_h^{(3)}=\begin{pmatrix}1&1&1\\-h&0&h\\h^2/2&0&h^2/2\end{pmatrix},
 \qquad (C_h^{(3)})^{-1}=
 \begin{pmatrix}0&-1/(2h)&1/h^2\\1&0&-2/h^2\\0&1/(2h)&1/h^2\end{pmatrix},
 \quad\det C_h^{(3)}=h^3.
\end{equation}
For later use allow arbitrary inverse strengths $b_-,b_c,b_+$.
Write $g=g_\lambda$, $t=G_\lambda(he)$, $v=G_\lambda(2he)$.
Direct multiplication gives
\begin{equation}\label{eq:full3}
 (C_h^{(3)})^{-T}Q_h^{(3)}(C_h^{(3)})^{-1}=
 \begin{pmatrix}
 b_c+g&0&\dfrac{2(t-g-b_c)}{h^2}\\[5pt]
 0&\dfrac{b_-+b_++2(g-v)}{4h^2}&\dfrac{b_+-b_-}{2h^3}\\[5pt]
 \dfrac{2(t-g-b_c)}{h^2}&\dfrac{b_+-b_-}{2h^3}&
 \dfrac{b_-+4b_c+b_++6g-8t+2v}{h^4}
 \end{pmatrix}.
\end{equation}
Here $Q_h^{(3)}=\diag(b_-,b_c,b_+)+(G_\lambda(x_i-x_j))_{i,j}$.
For $b_-=b_+=b_o$ the odd coordinate decouples. The even block is
\begin{align*}
 A_{h,11}&=b_c+g, & A_{h,12}&=2(t-g-b_c)/h^2,\\
 A_{h,22}&=(2b_o+4b_c+6g-8t+2v)/h^4,
 &k_{h,\lambda}&=(b_o+g-v)/(2h^2).
\end{align*}
Consequently the resolvent is exactly
\begin{equation}\label{eq:trifinite}
 R_h^{(3)}(\lambda)=R_0(\lambda)-U_{h,\lambda}A_h^{-1}U_{h,\bar\lambda}^*
       -D_{h,\lambda}^{s}k_{h,\lambda}^{-1}(D_{h,\bar\lambda}^{s})^*,
 \quad U_{h,\lambda}=(G_\lambda,E_{h,\lambda}).
\end{equation}

\emph{Even block and its logarithmic remainder.}
Substitute \eqref{eq:triscaling}. The first entry is $g+b_c$; the second is
$-2\beta$ minus the left side of \eqref{eq:renorm}. The last is
$\vartheta$ plus the fourth difference
\[
 \frac{6(G_\lambda-G_1)(0)-8(G_\lambda-G_1)(he)
                             +2(G_\lambda-G_1)(2he)}{h^4}.
\]
The quartic term in \eqref{eq:local} gives the fourth Fourier moment
$I_4(\lambda)$ in \eqref{eq:I4}. The next term is of order
$h^2(1+|\log h|)$ after division by $h^4$. Its logarithmic coefficient gives
\begin{equation}\label{eq:A22sharp}
 A_{h,22}-A_{\lambda,22}
 =\frac{5(\lambda-1)}{192\pi}h^2\log(1/h)+O(h^2).
\end{equation}
Indeed the weights $-8,2$ on $h^6,(2h)^6$ give $-8+128=120$,
and $120/4608=5/192$. The constant and quadratic terms cancel. Thus
\begin{equation}\label{eq:triexp}
 A_h=A_\lambda+O(h^2L_h),\qquad
 k_{h,\lambda}=\frac3{8\pi}\log(1/h)+O(1).
\end{equation}
To verify the second assertion, put $c_0=1/(8\pi)$ and $L=\log(1/h)$.
The local expansion yields
\[
 b_c=-c_0h^2L+O(h^2),\qquad b_o=2c_0h^2L+O(h^2),\qquad
 g-G_\lambda(2he)=4c_0h^2L+O(h^2).
\]
This proves the odd coefficient and the signs $b_o>0>b_c$.

For nonreal $\lambda$, subtraction of the Fourier multipliers at
$\lambda$ and $\bar\lambda$ gives
\[
 \operatorname{Im}A_\lambda=-(\operatorname{Im}\lambda)U_\lambda^*U_\lambda.
\]
For the off-diagonal entry this uses the negative second moment
$-I_2$, and for the lower-right entry the fourth moment $I_4$.
The resulting Gram matrix is positive definite: the multipliers
$1$ and $-\xi_e^2$ are independent. Hence $A_\lambda$ is invertible.
Equations \eqref{eq:trifinite}--\eqref{eq:triexp} prove convergence to
\eqref{eq:Eres}, which is the resolvent of a self-adjoint operator by
Lemma~\ref{lem:pseudo}.

\emph{Schur complement and common lower bound.}
Choose $\lambda_0$ so large that $A_{22}<0$. Since $A_{11}=g>0$,
$\det A<0$. For small $h$, $k_h>0$ and $\det A_h<0$.
The positive outer-site block is
\[
 Q_{++}=\begin{pmatrix}b_o+g&v\\v&b_o+g\end{pmatrix}>0,
 \qquad \det Q_{++}=2h^2k_h(b_o+g+v).
\]
Positivity follows from $b_o>0$ and positivity of the free Gram matrix.
Since $\det Q_h^{(3)}=h^6k_h\det A_h$, its center Schur complement is
\begin{equation}\label{eq:triS}
 S_h=\frac{h^4\det A_h}{2(b_o+g+v)},\qquad
 \frac{S_h}{h^4}\longrightarrow\frac{\det A_{\lambda_0}}{4g_{\lambda_0}}<0.
\end{equation}
In particular $Q_h^{(3)}$ has two positive and one negative eigenvalues.
In the positive energy $a_+$ obtained by adding the outer interactions,
let $r_h$ be the squared norm of center evaluation. The evaluation
representer calculation gives $b_c+r_h=S_h$, so
\[
 q_h^{(3)}[u]+\lambda_0\norm u_2^2
 \ge\frac{-S_h}{|b_c|}a_+[u]
 \ge c\frac{h^2}{L_h}\norm u_{H^2}^2\ge0.
\]
After increasing the shift by $1$, this proves the useful Schur estimate
\begin{equation}\label{eq:tricoerc}
 q_h^{(3)}[u]+\lambda\norm u_2^2
       \ge\norm u_2^2+c\frac{h^2}{L_h}\norm u_{H^2}^2.
\end{equation}

\emph{Optimal energy coercivity.} At a fixed sufficiently large shift,
use the relative form matrix
$F_h=I+M^{1/2}\diag(b_o^{-1},b_c^{-1},b_o^{-1})M^{1/2}$,
as justified in Theorem~\ref{thm:stable}. It commutes with reflection
of the two outer coordinates. Its odd eigenvalue is
$\nu_o=1+(g-v)/b_o\to3$. Its full determinant is
$h^6k_h\det A_h/(b_o^2b_c)$, so the product of its two even eigenvalues is
\[
 \delta_e\sim-\frac{\det A_\lambda}{4c_0^2L^2}.
\]
The exact prescription and the local expansion imply
$b_o+2b_c=-4c_0(\log2)h^2+o(h^2)$, and hence
\[
 2/b_o+1/b_c\sim\frac{2\log2}{c_0h^2L^2}.
\]
Subtracting the bounded odd eigenvalue from
$\tr F_h=3+g(2/b_o+1/b_c)$ gives the even trace
$t_e\sim2g\log2/(c_0h^2L^2)$. Since $\delta_e/t_e^2\to0$,
the smaller even eigenvalue is asymptotic to
$\delta_e/t_e=-\det A_\lambda\,h^2/(8c_0g\log2)$.
It is smaller than the other eigenvalues and than $1$, proving
\eqref{eq:optimaltri}. In particular a fixed larger shift gives
$q_h^{(3)}+\lambda\norm\cdot_2^2
\ge\norm\cdot_2^2+ch^2\norm\cdot_{H^2}^2$.
The extra logarithm in the weaker \eqref{eq:tricoerc} is not optimal.

\emph{Leading operator and its norm constant.}
The even contribution in \eqref{eq:trifinite} differs from its limit by
$O(h)$, by source convergence and \eqref{eq:triexp}. Its odd contribution
is $-D_h^sk_h^{-1}(D_h^s)^*$ at a real shift. The second estimate in
\eqref{eq:triexp} proves \eqref{eq:trirate}. Finally,
\begin{equation}\label{eq:Dnorm}
 \norm{D_\lambda}_2^2=\frac1{4\pi}\int_0^\infty
                   \frac{r^3}{(r^4+\lambda)^2}\dd r
                     =\frac1{16\pi\lambda}.
\end{equation}
Thus the limiting rank-one norm is $1/(6\lambda)$, and it is nonzero.

Since $U_{\bar\lambda}^*$ is onto and $A_\lambda$ is
invertible, some limiting resolvent range vector has a nonzero
$E_\lambda$ coefficient. Its other terms are in $H^1$, whereas
$E_\lambda\notin H^1$. This proves the domain assertion.
\end{proof}

For the smooth local triplet with the same tuning and
$\varepsilon=h^p$, $p>2$, Theorem~\ref{thm:hessianwidth} below gives
the same limit, common lower bound and leading operator: indeed its
width invariant satisfies $T_h=(p-1)h^{2p-4}\log(1/h)$, so
$T_h=o(1/\log(1/h))$. This is a consequence of the complete width
theorem, not an extra assumption in its proof.

\begin{proposition}[Three sites cannot activate all three collinear jets]
\label{prop:rankthree}
Allow arbitrary real diagonal strengths at $-he,0,he$, without imposing
symmetry or a common lower bound. A self-adjoint norm resolvent limit
cannot have rank-three correction in the channels
$G_\lambda,\partial_eG_\lambda,\partial_e^2G_\lambda$.
\end{proposition}
\begin{proof}
The three-column source map converges in norm to the injective map
$(G,D,E)$. Thus \eqref{eq:left}--\eqref{eq:recoverB} imply convergence
of the full transformed inverse. Rank three makes its limit invertible,
so the full matrix \eqref{eq:full3} itself must converge.
Put $c_0=1/(8\pi)$. Boundedness of its $(1,3)$ entry forces
$b_c=c_0h^2\log h+O(h^2)$. Boundedness of its $(2,2)$ entry forces
\[
 b_-+b_+=2(G_\lambda(2he)-g_\lambda)+O(h^2)
                      =8c_0h^2\log h+O(h^2).
\]
The numerator of its $(3,3)$ entry is consequently
\[
 b_-+4b_c+b_++6g_\lambda-8G_\lambda(he)+2G_\lambda(2he)
                       =12c_0h^2\log h+O(h^2).
\]
The kernel combination has no $h^2\log h$ term. Division by $h^4$
cannot converge, a contradiction. This proof includes nonsymmetric
strengths and uses no lower bound.
\end{proof}

\paragraph{Minimality.}
Three sites suffice by Theorem~\ref{thm:triplet}. Two scalar sites
cannot produce this Hessian interaction, for any directions or collision
scales: by Theorem~\ref{thm:derivativeclosure} a channel with the
quadratic symbol $\xi_e^2$ requires at least three sites. This does not
concern other approximation classes.

\subsection{A complete width classification in the separated-site regime}
Keep exactly the triplet strengths \eqref{eq:triscaling}, with fixed
$\beta,\vartheta$, and allow any positive width $\varepsilon(h)$ such
that $\varepsilon/h\to0$. Define
\begin{equation}\label{eq:widthinvariant}
 L=\log(1/h),\qquad \ell_h=\log(h/\varepsilon),\qquad
 T_h=\frac{\varepsilon^2}{h^4}\ell_h.
\end{equation}
This nonnegative quantity is an invariant of this declared regularization
scheme, not a coordinate-free invariant of arbitrary microscopic models.

\begin{theorem}[Sharp Hessian smoothing law]\label{thm:hessianwidth}
For all separated widths $\varepsilon/h\to0$, the smooth local triplets
with \eqref{eq:triscaling} have a common lower bound. They have a norm
resolvent limit if and only if $T_h$ has a limit in $[0,\infty]$.
The limit is
\begin{equation}\label{eq:hessianwidthlimit}
 \begin{cases}
 H^{(E)}_{\beta,\vartheta-3\kappa_\rho\tau},&T_h\to\tau<\infty,\\
 H_{\rm pin},&T_h\to\infty.
 \end{cases}
\end{equation}
In particular the originally prescribed Hessian limit is retained
exactly when $T_h\to0$.

If $T_h\to0$, at a fixed sufficiently large real shift, put
$P_2=\diag(0,1)$ and
$\mathcal E_\lambda=U_\lambda A_\lambda^{-1}P_2
A_\lambda^{-1}U_\lambda^*$. Then the local resolvent has the expansion
\begin{align}
 (\widetilde H_{h,\varepsilon}^{(3)}+\lambda)^{-1}
 -(H^{(E)}_{\beta,\vartheta}+\lambda)^{-1}
 ={}&-\frac{8\pi}{3L}D_\lambda D_\lambda^*
              -3\kappa_\rho T_h\mathcal E_\lambda
                         +o(L^{-1}+T_h).\label{eq:hessiantwoscale}
\end{align}
Consequently its norm equals
\begin{equation}\label{eq:hessianwidthnorm}
 \max\left\{\frac1{6\lambda L},\,
 3\kappa_\rho T_h\norm{U_\lambda A_\lambda^{-1}e_2}_2^2\right\}
                                 +o(L^{-1}+T_h),
 \qquad e_2=(0,1)^T.
\end{equation}
\end{theorem}
\begin{proof}
Use the averaged local traces of Section~\ref{sec:loc-smooth} and bars for their matrices.
Their exact entries are
\begin{align*}
 \overline A_{11}&=A_{h,11}+M_{\lambda,\varepsilon}(0)-g_\lambda,&
 \overline A_{12}&=A_{h,12}+I_1(h,\varepsilon),\\
 \overline A_{22}&=A_{h,22}-I_2(h,\varepsilon),&
 \overline k_h&=k_h-I_1(2h,\varepsilon).
\end{align*}
These follow by applying the same stencils to the two Gram kernels.
In particular no bound on an untransformed entry is substituted for its
actual divided difference. Lemma~\ref{lem:filter} gives
\begin{align}
 \overline A_{11}&=g_\lambda+o(1),&
 \overline A_{12}&=A_{\lambda,12}+o(1),\nonumber\\
 \overline A_{22}&=A_{\lambda,22}-3\kappa_\rho T_h
                       +O(T_h/\ell_h)+O(h^2L),&
 \overline k_h&=3c_0L+O(1)+o(1),\quad c_0=1/(8\pi).
                                                        \label{eq:averagedhessian}
\end{align}
Indeed $I_1(h,\varepsilon)=O(h^2T_h)=o(1)$, because
$h^2T_h=(\varepsilon/h)^2\log(h/\varepsilon)\to0$.
The source maps $P_\varepsilon(G,D_h^s,E_h)$ converge in norm to
$(G,D,E)$. These identities already give the averaged limits in
\eqref{eq:hessianwidthlimit}.

We next prove uniform lower bounds and a sufficiently accurate local
comparison for this whole range of widths. Choose a fixed shift so large
that $A_{\lambda,22}<-1$. Since $I_2\ge0$, for small $h$ the even
matrix has positive first entry and negative last entry, and
\[
                   -\det\overline A\asymp1+T_h.
\]
Its odd scalar is positive. The outer-site Schur argument used above
therefore proves a common lower bound for the averaged forms.

For the required quantitative coercivity, let $\overline F$ be their
relative form matrix. Write $m_0=M_{\lambda,\varepsilon}(0)$ and
$m_2=M_{\lambda,\varepsilon}(2he)$. Its odd eigenvalue is
\[
 \overline\nu_o=1+(m_0-m_2)/b_o\longrightarrow3.
\]
The even trace and determinant obey
\begin{align*}
 \overline t_e&=3+m_0(2/b_o+1/b_c)-\overline\nu_o
          \sim\frac{2g_\lambda\log2}{c_0h^2L^2},\\
 \overline\delta_e&=
 \frac{h^6\overline k_h\det\overline A}
                  {b_o^2b_c\overline\nu_o}
          \sim-\frac{\det\overline A}{4c_0^2L^2}.
\end{align*}
The scalar strength expansions in Theorem~\ref{thm:triplet} justify
the trace; the determinant follows from the source congruence.
Moreover $\overline\delta_e/\overline t_e^2\to0$, since it is
$O(h^4L^2(1+T_h))$ and $h^2T_h\to0$.
Thus the smaller even eigenvalue, and the optimal relative constant, are
\begin{equation}\label{eq:widthcoercivity}
 \overline\nu_h\sim
       -\frac{h^2\det\overline A}{8c_0g_\lambda\log2}
             \asymp h^2(1+T_h).
\end{equation}
It tends to zero, while the other even eigenvalue diverges, the odd
one tends to $3$, and the energy complement has relative value $1$.
Increase the fixed shift by one to give a separate lower bound by $I$.

The unsigned variance bound now has relative constant
\begin{equation}\label{eq:widthvariance}
 \eta_{h,\varepsilon}\le
 C\frac{\varepsilon^2L_\varepsilon}{h^4L(1+T_h)}
 =C\frac{T_h}{1+T_h}\frac{1+L+\ell_h}{L\ell_h}
 \le C(L^{-1}+\ell_h^{-1})\longrightarrow0.
\end{equation}
Here $\sum|\alpha_j|=O(h^{-2}/L)$ is essential.
Equation~\eqref{eq:formpert} proves a common lower bound and that the
local and averaged resolvents differ by $O(\eta_{h,\varepsilon})$.
This establishes both stated limits for the ordinary multiplication
potentials, not merely for the auxiliary averaged model.

To prove necessity, every sequence $h\downarrow0$ has a subsequence
on which $T_h$ converges in the compact extended half-line. Its
resolvent limit is one of \eqref{eq:hessianwidthlimit} by the proof just
given. Distinct finite values of $\tau$ give distinct inverse even
matrices and hence distinct operators, by the continuous left inverses
of $U_\lambda$. The pinned correction has rank one, whereas every
finite-$\tau$ correction has rank two. Thus a convergent operator
family has only one possible cluster value for $T_h$, which forces
convergence of $T_h$ itself.

For the last assertions suppose $T_h\to0$. Boundedness of $T_h$
forces $\ell_h\ge L/2$ eventually: otherwise
$T_h=e^{2L-2\ell_h}\ell_h\to\infty$ along that subsequence.
Hence \eqref{eq:widthvariance} improves to $O(T_h/L)$.
To bound the change of the source columns under convolution, the
second-difference multiplier gives
\[
 \|E_{h,\lambda}\|_{H^1}^2
 \le C_\lambda\int_0^\infty
 \frac{(1+r^2)\min(r^4,h^{-4})}{(1+r^4)^2}\,r\,\dd r
 \le C_\lambda(1+L).
\]
The three ranges $r<1$, $1<r<1/h$ and $r>1/h$ contribute
$O(1)$, $O(L)$ and $O(1)$, respectively. The $G$ and $D_h^s$
columns are uniformly bounded in $H^1$. The translation estimate
for $P_\varepsilon-I$ therefore bounds the change of the finite-$h$
source map by $O(\varepsilon\sqrt{1+L})=O(h^2)$, using
$\varepsilon^2=h^4T_h/\ell_h$.
Equations \eqref{eq:averagedhessian} therefore give the even-matrix
perturbation $-3\kappa_\rho T_hP_2+o(T_h)+O(h^2L)$.
The first-order inverse formula gives its contribution
$-3\kappa_\rho T_h\mathcal E_\lambda$ to the resolvent difference.
The point-source error is $O(h)=o(1/L)$; the odd block gives the first
term of \eqref{eq:hessiantwoscale}. The quadratic inverse remainders
are $O((T_h+h^2L)^2)=o(T_h+L^{-1})$, including when $T_h$ is
smaller than every fixed power of $h$. The local variance error is
also $o(T_h+L^{-1})$.
This proves the expansion. The odd vector $D_\lambda$ is orthogonal
to the even vector $U_\lambda A_\lambda^{-1}e_2$ by reflection parity.
The two leading negative rank-one terms act on orthogonal subspaces.
Their norm is their larger rank-one norm, proving
\eqref{eq:hessianwidthnorm} and its nonzero constants.
\end{proof}

\begin{corollary}[The two sharp logarithmic width thresholds]
\label{cor:widththresholds}
Within $\varepsilon/h\to0$, the original Hessian limit is retained
if and only if
\[
                  \varepsilon=o\bigl(h^2/\sqrt{\log(1/h)}\bigr).
\]
For $\varepsilon=s h^2/\sqrt{\log(1/h)}$, $s>0$, the limit instead has
$\vartheta$ replaced by $\vartheta-3\kappa_\rho s^2$.
For powers $\varepsilon=h^p$, $p>1$, the original limit occurs
exactly when $p>2$; $1<p\le2$ gives $H_{\rm pin}$.
Among widths retaining the original limit, its leading operator
\eqref{eq:trirate} is retained if and only if
\[
                         \varepsilon=o(h^2/\log(1/h)).
\]
At $\varepsilon=s h^2/\log(1/h)$ the additional leading operator is
$-3\kappa_\rho s^2\mathcal E_\lambda$ after multiplication by
$\log(1/h)$.
\end{corollary}
\begin{proof}
For bounded $T_h$, the preceding proof gives $\ell_h\ge L/2$.
Thus $T_h\to0$ implies $\varepsilon\sqrt L/h^2\to0$.
Conversely write $\varepsilon=h^2b_h/\sqrt L$, $b_h\to0$; then
\[
 T_h=b_h^2\left[1+\frac{\log L}{2L}+\frac{\log(1/b_h)}L\right]\to0,
\]
because $b^2\log(1/b)\to0$. The critical substitution gives $T_h\to s^2$.
For powers, $T_h=(p-1)h^{2p-4}L$.
The expansion and parity in Theorem~\ref{thm:hessianwidth} show that
the original leading operator is retained exactly when $LT_h\to0$.
The same argument with $\varepsilon=h^2b_h/L$ makes this equivalent
to $b_h\to0$. At the stated critical width $LT_h\to s^2$.
\end{proof}

\begin{proposition}[Spectrum of the directional Hessian interaction]
\label{prop:Espectrum}
$H^{(E)}_{\beta,\vartheta}$ has essential spectrum $[0,\infty)$.
It has one simple negative eigenvalue if and only if
$\vartheta+3/64>0$; otherwise it has no negative eigenvalues.
The eigenvalue $-\lambda_*$ is determined by
\begin{equation}\label{eq:Eeigen}
 \vartheta+\frac3{64}
 =\sqrt{\lambda_*}\left[\frac3{64}
       +8\left(-2\beta+\frac{\log\lambda_*}{16\pi}\right)^2\right].
\end{equation}
Neither $H^{(E)}_{\beta,\vartheta}$ nor $H_{a,\theta}$ has a zero or
positive eigenvalue. For $\beta=\vartheta=0$, the triplet point models
have the exact eigenvalue $-1$ for every sufficiently small $h$, with
eigenfunction $E_{h,1}$; the limiting eigenfunction is $E_1$.
\end{proposition}
\begin{proof}
The scalar even Schur complement is
\[
 S(\lambda)=\vartheta+\frac3{64}
   -\sqrt\lambda\left[\frac3{64}+8q(\lambda)^2\right],\qquad
 q(\lambda)=-2\beta+\frac{\log\lambda}{16\pi}.
\]
The derivative of the subtracted term is
\[
 \frac1{2\sqrt\lambda}
 \left[\frac3{64}+8q^2+\frac{2q}\pi\right]>0,
\]
because the bracket is
$8(q+1/(8\pi))^2+3/64-1/(8\pi^2)>0$.
That term tends to zero at $\lambda=0$ and to infinity at infinity.
This proves the exact existence condition and uniqueness.
At a root $A_{11}=g>0$, so the kernel is one dimensional; the
strict derivative gives a simple rank-one resolvent pole, as in the
proof of Proposition~\ref{prop:pairnegative}. Its eigenfunction is a nonzero multiple of
$U_\lambda(-A_{12}/g,1)^T$.
Weyl's theorem \cite[Section~6.4]{Teschl} gives
$\sigma_{\rm ess}=[0,\infty)$, and the explicit even resolvent excludes
any other negative spectral point.

For an eigenfunction at energy $E\ge0$, the explicit resolvent formula
gives distributionally
$(\Delta^2-E)u=w_0\delta_0+w_1\partial_e\delta_0$ in the pair
case, or $w_0\delta_0+w_2\partial_e^2\delta_0$ in the Hessian case.
For $E>0$, its Fourier numerator must vanish on the whole circle
$|\xi|^4=E$ to give an $L^2$ quotient; if it does not vanish at one
point, integration across a small circle arc diverges as
$\int|r-E^{1/4}|^{-2}\dd r$. Neither a nonzero linear polynomial
$w_0+iw_1\xi_e$ nor $w_0-w_2\xi_e^2$ vanishes on that whole circle.
If the numerator is zero, an $L^2$ function supported on the circle
is zero. At $E=0$, the lowest nonzero polynomial degree is at most
two, and its quotient by $|\xi|^4$ fails $L^2$ at the origin.
Thus there are no such eigenfunctions.

Finally at $\beta=\vartheta=0$ and $\lambda=1$, the exact even block
has $A_{h,12}=A_{h,22}=0$. Its kernel vector is $(0,1)$ and the
associated source is $E_{h,1}$. The point-model eigenvalue correspondence
in Section~\ref{sec:loc-model} gives the exact eigenvalue. Its limiting version follows
from \eqref{eq:Eeigen}, or directly from $\ker A_1=\operatorname{span}\{(0,1)^T\}$.
\end{proof}

\section{Spectral and geometric obstructions}\label{sec:loc-obstructions}

\subsection{A stability criterion for scalar configurations}\label{sec:loc-inertia}

For this section let $x_1,\ldots,x_N$ be distinct sites, let
$B=\diag(b_1,\ldots,b_N)$ with real $b_j\ne0$, and let
$Q(\lambda)=B+M_\lambda$, $\lambda>0$, as in
Section~\ref{sec:loc-model}.
Write $n_+(T),n_-(T),n_0(T)$ for the positive, negative and zero inertia
indices of a finite Hermitian matrix. All negative-eigenvalue counts below
include multiplicity.

\begin{theorem}[Inertia and a common lower bound]\label{thm:inertia}
For the form operator $H$ with these sites,
\begin{equation}\label{eq:inertiaexact}
 N_H(-\infty,-\lambda)
   =n_-(B)-n_-(Q(\lambda))-n_0(Q(\lambda)).
\end{equation}
Consequently a family of such finite scalar models has a common lower
bound for small $h$ if and only if some fixed $\lambda_0>0$ has a zero
right side in \eqref{eq:inertiaexact} for every sufficiently small $h$.
Equivalently, after increasing the shift if necessary, $Q_h(\lambda_0)$
is invertible and $n_-(Q_h(\lambda_0))=n_-(B_h)$ eventually.
\end{theorem}
\begin{proof}
In the free energy space the form on the evaluation representers is
$F=I+M^{1/2}B^{-1}M^{1/2}$, $M=M_\lambda$, and is the identity on
their orthogonal complement, as in Theorem~\ref{thm:stable}.
The min--max principle \cite[Section~4.4]{Teschl} gives
$N_H(-\infty,-\lambda)=n_-(F)$. Taking the two Schur complements
of $\left(\begin{smallmatrix}B&M^{1/2}\\ M^{1/2}&-I\end{smallmatrix}\right)$
in Haynsworth's inertia formula \cite{Hay} gives
$n_+(Q)=n_+(B)+n_-(F)$ and $n_0(Q)=n_0(F)$.
Counting the positive eigenvalues of $Q$ proves \eqref{eq:inertiaexact},
including its nullity term. A zero count is equivalent to
$H\ge-\lambda$ and a strictly larger shift removes the kernel.
\end{proof}

\begin{proposition}[Positive-site Schur complement]\label{prop:generalSchur}
Partition the sites into $b_j>0$ and $b_j<0$. Put
\[
 Q_{++}=B_++M_{++}>0,\quad
 R_-=M_{--}-M_{-+}Q_{++}^{-1}M_{+-},\quad
 S=B_-+R_-.
\]
With the empty positive block interpreted by deletion, the shifted form
is nonnegative exactly when $S\le0$. If
$S\le-\eta|B_-|$ for $0<\eta\le1$, then
\begin{equation}\label{eq:generalSchurbound}
 q[u]+\lambda\norm u_2^2\ge
 \eta\left(a_{0,\lambda}[u]+\sum_{j\in+}b_j^{-1}|u(x_j)|^2\right).
\end{equation}
If there are no negative sites, the shifted form is positive without a
Schur test.
\end{proposition}
\begin{proof}
After the positive sites have been added, let $a_+$ denote their energy
and let $\tau_-$ be negative-site evaluation. Subtracting the positive-site
resolvent correction gives its Gram matrix $R_-$. It is positive definite
because independent evaluations remain independent in the equivalent
fixed-configuration energy norm. The negative contribution is
$-\norm{|B_-|^{-1/2}\tau_-u}^2$.
Its squared operator norm relative to $a_+$ is the largest eigenvalue
of $|B_-|^{-1/2}R_-|B_-|^{-1/2}$. Therefore the full form is nonnegative
if and only if $R_-\le|B_-|$, or $S\le0$. The stronger hypothesis bounds
this squared norm by $1-\eta$, proving \eqref{eq:generalSchurbound}.
\end{proof}

\paragraph{Escaping negative eigenvalues.}
Let $I_h(\lambda)$ be the right side of \eqref{eq:inertiaexact}.
There are eigenvalues tending to $-\infty$ along some $h_n\downarrow0$
if and only if, for every $\Lambda>0$, $I_h(\Lambda)>0$ for arbitrarily small $h$.
Indeed, the latter condition permits a choice $h_n<1/n$ with an
eigenvalue below $-n$; the converse is immediate. This is a spectral
criterion for the declared scalar class, not a classification of the
asymptotic scalar strengths. It is invariant under any invertible source
congruence, since that congruence preserves the inertia of $Q_h$.

\subsection{Finite-site geometry and dihedral instability}\label{sec:loc-geometry}

We first isolate a source-space obstruction, without any stability
assumption. A Hessian direction means a homogeneous quadratic channel
$\sum_{i,j}T_{ij}\partial_i\partial_jG_\lambda$ modulo the monopole and
gradient channels, with real symmetric $T$.

\begin{theorem}[At least five fixed scalar sites are necessary]
\label{thm:fourlower}
Let $x_j(h)=ha_j$ with at most four distinct fixed reference points
$a_j\in\R^2$. If the scalar point models have a self-adjoint norm
resolvent limit, its correction cannot contain two independent Hessian
directions. In particular, two independent Hessians require at least five
sites in this fixed-configuration class. No sufficiency assertion for five
sites is made.
\end{theorem}
\begin{proof}
This is the case $N\le4$ of Theorem~\ref{thm:derivativeclosure}, which
holds for arbitrary site geometries and whose proof does not use this
section.
\end{proof}

For the cross order $(he_1,-he_1,he_2,-he_2)$ define, with translates
abbreviated by their site subscripts,
\[
 W_h=\left(\frac{G_{+x}+G_{-x}+G_{+y}+G_{-y}}4,
 \frac{G_{+x}-G_{-x}}{2h},\frac{G_{+y}-G_{-y}}{2h},
 \frac{G_{+x}+G_{-x}-G_{+y}-G_{-y}}{h^2}\right).
\]
Its limit is $(G,D_1,D_2,\partial_1^2G-\partial_2^2G)$, with errors
$O(h^2),O(h),O(h),O(h)$ in $L^2$, as in the proof of
Theorem~\ref{thm:triplet}.
The exact source matrices are
\begin{equation}\label{eq:crossC}
 C_h=\begin{pmatrix}1&1&1&1\\h&-h&0&0\\0&0&h&-h\\
 h^2/4&h^2/4&-h^2/4&-h^2/4\end{pmatrix},\qquad
 C_h^{-1}=\begin{pmatrix}
 1/4&1/(2h)&0&1/h^2\\1/4&-1/(2h)&0&1/h^2\\
 1/4&0&1/(2h)&-1/h^2\\1/4&0&-1/(2h)&-1/h^2
 \end{pmatrix}.
\end{equation}
Put $g=G_\lambda(0)$, $t=G_\lambda(2he_1)$,
$s=G_\lambda(h(e_1-e_2))$, $S_x=b_1+b_2$, $S_y=b_3+b_4$,
$T_x=b_1-b_2$, $T_y=b_3-b_4$, and $\mu_0=g+t+2s$,
$\mu_H=g+t-2s$. Write $K_h^\times=C_h^{-T}Q_hC_h^{-1}$.
Direct multiplication gives the full congruence
\begingroup\setlength{\arraycolsep}{3.5pt}
\begin{equation}\label{eq:crossK}
 K_h^\times=
 \begin{pmatrix}
 \mu_0/4+(S_x+S_y)/16&T_x/(8h)&T_y/(8h)&(S_x-S_y)/(4h^2)\\[3pt]
 T_x/(8h)&[S_x+2(g-t)]/(4h^2)&0&T_x/(2h^3)\\[3pt]
 T_y/(8h)&0&[S_y+2(g-t)]/(4h^2)&-T_y/(2h^3)\\[3pt]
 (S_x-S_y)/(4h^2)&T_x/(2h^3)&-T_y/(2h^3)&
                   [S_x+S_y+4\mu_H]/h^4
 \end{pmatrix}.
\end{equation}
\endgroup
The source coordinates transform as $A_1\oplus E\oplus B_1$ under $D_4$:
the average is invariant, $(D_1,D_2)$ transform as the coordinate vector,
and the last column changes sign under a quarter-turn and is unchanged
by reflection in a coordinate axis. These transformation rules explicitly
specify the representation labels. There is no independent trace-Hessian
coordinate, which shows the cross obstruction directly.

\begin{proposition}[Full dihedral symmetry forces escaping states]
\label{prop:D4escape}
If $b_1=b_2=b_3=b_4=b(h)$ and a self-adjoint norm resolvent limit has a
nonzero Hessian component, then at least three negative eigenvalues escape
to $-\infty$ as $h\downarrow0$. In particular such a family cannot have a
common lower bound.
\end{proposition}
\begin{proof}
In \eqref{eq:crossK} the diagonal entries become
\[
 (b+\mu_0)/4,\quad(b+g-t)/(2h^2),\quad(b+g-t)/(2h^2),
                    \quad4(b+\mu_H)/h^4.
\]
At a fixed nonreal parameter, coefficient recovery and a nonzero Hessian
component force $h^4/[4(b+\mu_H)]$ to converge to a nonzero number.
The local expansion \eqref{eq:local} gives
\[
 \mu_H=2c_0(\log2)h^2-\frac{\sqrt\lambda}{64}h^4+O(h^6L_h),
 \qquad c_0=\frac1{8\pi}.
\]
Thus $b(h)=-2c_0(\log2)h^2+O(h^4)<0$.
At any fixed real $\Lambda>0$, $b+\mu_0\to4g_\Lambda>0$, while
$b+g_\Lambda-G_\Lambda(2he_1)=4c_0h^2\log(1/h)+O(h^2)>0$.
Hence $Q_h(\Lambda)$ has at least three positive eigenvalues, and
$n_-(B_h)=4$, $n_-(Q_h(\Lambda))+n_0(Q_h(\Lambda))\le1$.
Theorem~\ref{thm:inertia} gives at least three eigenvalues below $-\Lambda$
for every sufficiently small $h$. Since $\Lambda$ is arbitrary, the three
lowest negative eigenvalues tend to $-\infty$.
\end{proof}

\paragraph{Accessible sources and scalar realizations.}
A mixed derivative is also an $L^2$ divided-difference source: for corner
sites $(\pm h,\pm h)$,
\[
 \frac{G_{++}-G_{+-}-G_{-+}+G_{--}}{4h^2}
           \longrightarrow\partial_1\partial_2G\quad\hbox{in }L^2.
\]
The Fourier multiplier is
$-\sin(h\xi_1)\sin(h\xi_2)/h^2$, bounded by $|\xi_1\xi_2|$,
so dominated convergence applies after division by $|\xi|^4+\lambda$.
This proves source convergence only; no stable scalar-coupling construction
for a mixed-derivative limit is asserted here.

\subsubsection{Adding the central site}
For the five sites $0,\pm he_1,\pm he_2$, define the central divided
differences $D_{x,h},D_{y,h},E_{x,h},E_{y,h}$ as in Section~\ref{sec:loc-hessian}, in each
coordinate direction. In the site order $0,+x,-x,+y,-y$, the source
map is $\Phi_h=V_hC_h^{(5)}$, where
\[
 V_h=(G,D_{x,h},D_{y,h},E_{x,h},E_{y,h}),\qquad
 C_h^{(5)}=\begin{pmatrix}
 1&1&1&1&1\\0&h&-h&0&0\\0&0&0&h&-h\\
 0&h^2/2&h^2/2&0&0\\0&0&0&h^2/2&h^2/2
 \end{pmatrix}.
\]
Its inverse columns are $e_0$, $(0,1,-1,0,0)^T/(2h)$,
$(0,0,0,1,-1)^T/(2h)$, $(-2,1,1,0,0)^T/h^2$ and
$(-2,0,0,1,1)^T/h^2$.
These identities follow by recovering the five translates from the
definitions of the four differences. The five limiting source columns
are independent, and two Hessian directions are available at the level
of source span.

Give all four outer sites the same inverse strength $b_o$ and the
centre the inverse strength $b_c$. Put $u=G_\lambda(h)$,
$t=G_\lambda(2h)$, $s=G_\lambda(\sqrt2h)$ and $g=g_\lambda$.
In the source coordinates
\[
 G,\quad (E_{x,h}+E_{y,h})/\sqrt2,\quad
 D_{x,h},D_{y,h},\quad (E_{x,h}-E_{y,h})/\sqrt2,
\]
the congruence splits into the following even block and scalars:
\begin{align}
 K_{A_1}&=\begin{pmatrix}
 b_c+g&2\sqrt2(u-g-b_c)/h^2\\
 2\sqrt2(u-g-b_c)/h^2&
 [8b_c+2b_o+10g-16u+2t+4s]/h^4
 \end{pmatrix},\nonumber\\
 k_E&=(b_o+g-t)/(2h^2)\quad\hbox{(twice)},\qquad
 k_{B_1}=2(b_o+g+t-2s)/h^4.                 \label{eq:starblocks}
\end{align}
For verification, the squared second-difference free numerator is
$6g-8u+2t$, and the cross numerator of the two Hessians is
$4g-8u+4s$. Their sum and difference give the lower-right entry of
$K_{A_1}$ and $k_{B_1}$; the centre weights are $-2$ in each stencil.
Odd free cross terms cancel, and their squared numerator is $2(g-t)$.
These computations prove every entry of \eqref{eq:starblocks}.
The representation is $2A_1\oplus E\oplus B_1$: the two invariant
sources are $G$ and the trace Hessian, the dipoles form $E$, and the
anisotropic Hessian changes sign under a quarter-turn.

\begin{proposition}[The symmetric star still has escaping states]
\label{prop:starescape}
Suppose a self-adjoint norm resolvent limit of these fully symmetric
five-site models has a nonzero anisotropic Hessian component. Then its
two lowest eigenvalues tend to $-\infty$. This conclusion is independent
of the centre strength. In particular a common lower bound rules out
activation of both independent Hessians in this symmetric star.
\end{proposition}
\begin{proof}
At a fixed nonreal shift, the five source columns converge to an
injective map. The continuous left inverses of Section~\ref{sec:loc-converse} therefore
recover a convergent transformed inverse. Symmetry decouples its
$B_1$ scalar, so an anisotropic component forces $k_{B_1}^{-1}$ to
converge to a nonzero value. The expansion already used for the cross
gives
\[
 g+t-2s=2c_0(\log2)h^2-\sqrt\lambda h^4/64+O(h^6L_h).
\]
Consequently $b_o=-2c_0(\log2)h^2+O(h^4)<0$.

At any fixed real shift $\Lambda>0$, the two odd site vectors
$(0,1,-1,0,0)^T$ and $(0,0,0,1,-1)^T$ are simultaneous eigenvectors
of $B$ and $M_\Lambda$, with eigenvalues $b_o$ and
$g_\Lambda-G_\Lambda(2h)$, respectively. Each is thus an eigenvector
of the relative form matrix $I+M_\Lambda^{1/2}B^{-1}M_\Lambda^{1/2}$
with eigenvalue
\[
 1+\frac{g_\Lambda-G_\Lambda(2h)}{b_o}
              =-\frac{2}{\log2}\log(1/h)+O(1)<0.
\]
The shifted form has at least two negative directions. The finite-rank
index argument in Theorem~\ref{thm:inertia} gives at least two
eigenvalues below $-\Lambda$ for every sufficiently small $h$.
Since $\Lambda$ is arbitrary, the two lowest eigenvalues escape.
The centre does not enter either odd vector. The same odd argument
also applies if a zero centre strength is handled by deleting that site.
\end{proof}

\subsection{An attractive-site obstruction}\label{sec:loc-signcount}

\label{sec:signcount}
The number of attractive sites imposes an obstruction independent of the
moment geometry. We use the complete six-jet coordinates of
Theorem~\ref{fd:thm:planarM}: a regular jet is written $(a,d,H)$,
where $a\in\C$, $d\in\C^2$, and $H=(H_{11},H_{12},H_{22})\in\C^3$.
Let $E\subset\C^6$ be the active subspace of a self-adjoint relation and
let $\Theta=\Theta^*$ act on $E$. Its derivative rank is
\[
 r_{\rm der}(E)=\dim P_{\rm der}E,\qquad
 P_{\rm der}(a,d,H)=(d,H).
\]
In particular, two independent active Hessian coordinates give
$r_{\rm der}(E)\ge2$. This rank concerns the active relation, not the
number of nonzero entries of one matrix chart.

\begin{theorem}[Attractive sites are necessary for derivative channels]
\label{thm:attractiverank}
Suppose finite scalar point models $H_h$ have a common lower bound,
converge in norm resolvent sense to a self-adjoint extension $H$ of the
complete six-jet restriction at the origin, and have at most $p$
negative scalar strengths. If $E$ is the active subspace of $H$, then
\[
                         r_{\rm der}(E)\le p.
\]
Consequently a family with at most one attractive site cannot activate
two independent Hessian channels, irrespective of its number of
repulsive sites or its geometry. A full six-channel active relation
requires at least five attractive sites. The conclusion also holds for
smooth local realizations that are norm resolvent close to these point
models at one common real shift.
\end{theorem}
\begin{proof}
Choose a real $\lambda$ above the common lower bound in absolute value.
Theorem~\ref{thm:inertia} gives
$n_-(Q_h(\lambda))=n_-(B_h)\le p$. The source columns at distinct
sites are linearly independent: applying $\Delta^2+\lambda$ to a
linear combination gives independent point masses. The finite-rank
formula therefore gives
\[
 n_+\big((H_h+\lambda)^{-1}-R_0(\lambda)\big)
                  =n_-(Q_h(\lambda))\le p.                 \tag{*}
\]
Here $n_+$ counts strictly positive eigenvalues of a finite-rank
self-adjoint operator.

For completeness, the high-energy argument on the limiting active
space is explicit. Put
\[
 g_\lambda=(8\sqrt\lambda)^{-1},\quad
 L_\lambda=\frac{\log\lambda}{16\pi},\quad
 T_\lambda=\frac{1-\sqrt\lambda}{64}.
\]
The complete Weyl matrix in Theorem~\ref{fd:thm:planarM} has quadratic
form
\begin{align*}
 (v,M_{-\lambda}v)
 ={}&(g_\lambda-g_1)|a|^2-L_\lambda|d|^2
       +2L_\lambda\Re\big(\bar a(H_{11}+H_{22})\big)\\
 &+T_\lambda\big(3|H_{11}|^2+|H_{12}|^2+3|H_{22}|^2
                         +2\Re(\bar H_{11}H_{22})\big).
\end{align*}
The last parenthesis is positive definite. Choose a complement
$F\subset E$ of $E\cap\ker P_{\rm der}$, so $\dim F=r_{\rm der}(E)$
and $|d|^2+|H|^2\ge c_F|v|^2$ on $F$. Absorb the mixed term by
\[
 C\log\lambda\,|a|\,|H|
 \le c\sqrt\lambda\,|H|^2+
          C'\frac{(\log\lambda)^2}{\sqrt\lambda}|a|^2,
\]
choosing $c$ smaller than half the negative Hessian coefficient.
The bounded matrix $\Theta$ and the bounded monopole coefficient are
then dominated by the remaining negative derivative terms. Thus
$A_\lambda=P_EM_{-\lambda}|_E-\Theta$ is strictly negative on $F$
for every sufficiently large $\lambda$, and
$n_-(A_\lambda)\ge r_{\rm der}(E)$.

The relation resolvent formula is
\[
 (H+\lambda)^{-1}-R_0(\lambda)
       =-\gamma_E(-\lambda)A_\lambda^{-1}\gamma_E(-\lambda)^* .
\]
The map $\gamma_E$ is injective, so this difference has at least
$r_{\rm der}(E)$ strictly positive eigenvalues. Such a finite positive
index persists under sufficiently small operator-norm perturbations:
on any finite-dimensional positive spectral subspace the quadratic form
has a positive minimum on its unit sphere. Pass to the norm resolvent
limit in (*) to obtain the bound. The same norm-limit argument proves
the last assertion.
\end{proof}

\begin{proposition}[Derivative limits require loss of uniform bending coercivity]
\label{prop:bendingcoercivity}
Let $H_\varepsilon$ be self-adjoint operators on $L^2(\R^2)$ with
closed forms $q_\varepsilon$, whose form domains are contained in $H^2$.
Suppose that constants $c>0$ and $C\ge0$, independent of $\varepsilon$,
satisfy
\[
 q_\varepsilon[u]\ge c\|\Delta u\|_2^2-C\|u\|_2^2.
\]
If $H_\varepsilon$ converges in norm resolvent sense to a six-jet
point extension $H$, then its active subspace has derivative rank zero.
In particular, nonnegative scalar potentials cannot produce an active
dipole or Hessian channel. A stable derivative-channel realization must
lose every such uniform positive bending coefficient.
\end{proposition}
\begin{proof}
For $s>C$, the shifted form identity bounds $u_\varepsilon=
(H_\varepsilon+s)^{-1}f$ uniformly in $H^2$. Weak compactness in $H^2$
and norm resolvent convergence in $L^2$ give $\Dom(H)\subset H^2$.
In the direct-sum coordinates of the complete point restriction, a
domain element is a regular $H^4$ function plus Green derivatives
through order two. A combination with a nonzero derivative charge
cannot be in $H^2$. Indeed its Fourier transform has the form
$P(\xi)/(|\xi|^4+1)$ with a polynomial $P$ of degree one or two.
On a cone where its leading homogeneous polynomial is nonzero,
the squared $H^2$ integral at infinity bounds from below a multiple
of $\int_1^\infty r^{2\deg P-3}\,dr$, which diverges.
Adding an $H^4$ function cannot repair this failure. Every active
charge is attained by a domain element in the relation calculus;
therefore all derivative charges must be zero. Nonnegative potentials
satisfy the hypothesis with $c=1$ and $C=0$.
\end{proof}

\section{Stable multipoles and minimum site counts}\label{sec:loc-multipoles}

\subsection{Simultaneous Hessians from three radial triplets}\label{sec:loc-twohessians}

\label{sec:twohessians}
The fixed five-site question can be bypassed by using two collision
scales. We construct nine scalar sites, with three attractive centres
and six repulsive outer sites. The inner triplet scale is chosen after
the outer triangle scale. No power law relating these scales is asserted.

Put
\[
 a_0=(1,0),\quad a_1=(-1/2,\sqrt3/2),\quad
 a_2=(-1/2,-\sqrt3/2),\qquad p_j=s a_j,
 \qquad r=\sqrt3s.
\]
Each triplet is directed along its radial unit vector $a_j$.
Write $f_\lambda(r)=G_\lambda(r e_1)$ and define the two cross entries
\begin{align}
 m_s(\lambda)&=\tfrac34 f_\lambda''(r)
                       +\tfrac1{4r}f_\lambda'(r),\nonumber\\
 n_s(\lambda)&=\tfrac1{16}\left(9f_\lambda^{(4)}(r)
       -\frac6r f_\lambda'''(r)+\frac{15}{r^2}f_\lambda''(r)
       -\frac{15}{r^3}f_\lambda'(r)\right).\label{eq:radialcross}
\end{align}
They are respectively $\partial_{a_j}^2G_\lambda(p_i-p_j)$ and
$\partial_{a_i}^2\partial_{a_j}^2G_\lambda(p_i-p_j)$ for any $i\ne j$.

\begin{lemma}[Radial cancellation on an equilateral triangle]
\label{lem:radialcancel}
For each fixed $\lambda$ in the cut plane, as $s\downarrow0$,
\begin{align}
 m_s(\lambda)&=m_s(1)+\frac{\log\lambda}{16\pi}+O(s^2),\nonumber\\
 m_s(1)&=2c_0\log r+c_0(2\ell_1+1/2)+O(s^2),\nonumber\\
 n_s(\lambda)&=-\frac{3\sqrt\lambda}{128}
                         +O(s^2(1+|\log s|)),\label{eq:radialexp}\\
 g_\lambda-f_\lambda(r)&\sim3c_0s^2\log(1/s),
 \qquad c_0=\frac1{8\pi}.\nonumber
\end{align}
The last asymptotic is a complex asymptotic with positive real leading
coefficient; it holds also at fixed nonreal $\lambda$.
\end{lemma}
\begin{proof}
For $i\ne j$, let $n=(p_i-p_j)/r$, $e=a_i$ and $v=a_j$.
Then $e\cdot v=-1/2$, $(e\cdot n)^2=(v\cdot n)^2=3/4$, and
$(e\cdot n)(v\cdot n)=-3/4$. Twice differentiating a radial function
gives the first formula in \eqref{eq:radialcross}. Four differentiations
give
\[
 \partial_e^2\partial_v^2f
 =A(e\cdot n)^2(v\cdot n)^2
  +B\big((e\cdot n)^2+(v\cdot n)^2
                +4(e\cdot v)(e\cdot n)(v\cdot n)\big)
  +C\big(1+2(e\cdot v)^2\big),
\]
where $A=f^{(4)}-6f'''/r+15f''/r^2-15f'/r^3$,
$B=f'''/r-3f''/r^2+3f'/r^3$ and $C=f''/r^2-f'/r^3$.
Substitution gives the second formula in \eqref{eq:radialcross}.
Its differential expression annihilates $r^2\log r$ and $r^2$,
and takes the value $12$ on $r^4$. Thus \eqref{eq:local}, whose
convergent Bessel expansion may be differentiated on $r>0$, gives
$12(-\sqrt\lambda/512)=-3\sqrt\lambda/128$ with the stated remainder.
The first cross entry on $c_0r^2(\log r+\ell_\lambda-1)$ is
$c_0(2\log r+2\ell_\lambda+1/2)$. Since
$2c_0(\ell_\lambda-\ell_1)=\log\lambda/(16\pi)$, the remaining
assertions follow from the same expansion.
\end{proof}

\begin{lemma}[The inner limit of three separated triplets]
\label{lem:separatedtriplets}
Fix $s>0$ and real $\beta,\vartheta$. At each $p_j$ place the
triplet $p_j-ha_j,p_j,p_j+ha_j$, with the strengths
\eqref{eq:triscaling} using this common $\beta,\vartheta$.
As $h\downarrow0$, their resolvents converge in norm, at every nonreal
shift, to
\begin{equation}\label{eq:threeclusterres}
 R_s(\lambda)=R_0(\lambda)
          -\Gamma_{s,\lambda}K_s(\lambda)^{-1}
                                  \Gamma_{s,\bar\lambda}^*,
 \quad
 \Gamma_{s,\lambda}=(G_0,G_1,G_2,E_0,E_1,E_2),
\end{equation}
where $G_j=G_\lambda(\cdot-p_j)$ and
$E_j=\partial_{a_j}^2G_\lambda(\cdot-p_j)$. In this order the
diagonal and off-diagonal entries of the three blocks of $K_s$ are
\[
 \begin{array}{c|cc}
 &\text{diagonal}&\text{off-diagonal}\\\hline
 GG&g_\lambda&f_\lambda(r)\\
 GE&-2\beta+\log\lambda/(16\pi)&m_s(\lambda)\\
 EE&\vartheta+3(1-\sqrt\lambda)/64&n_s(\lambda).
 \end{array}
\]
If $K_s(\lambda_0)$ has exactly three negative eigenvalues and is
invertible at a real $\lambda_0>0$, then all sufficiently small inner
scales have $H_{s,h}\ge-\lambda_0$. Their resolvents at this shift also
converge to \eqref{eq:threeclusterres}.
\end{lemma}
\begin{proof}
Apply the exact three-site change of sources \eqref{eq:C3} separately
at each centre. The three odd diagonal entries are
$3\log(1/h)/(8\pi)+O(1)$, and the three even two-by-two diagonal
blocks, one for each centre, tend to the matrix $A_\lambda$ in
\eqref{eq:EA}. For different centres the
kernel is smooth in a fixed neighbourhood of their difference, so all
cross divided differences converge to the corresponding derivatives.
In particular, every off-diagonal entry involving an odd source remains
bounded for this fixed $s$. The odd block is therefore
$\frac{3\log(1/h)}{8\pi}I+O(1)$. At a nonreal shift it is complex
symmetric rather than Hermitian, but its real part is positive definite
for small $h$; since $\Re(x^*Yx)=x^*(\Re Y)x$ for a complex
symmetric $Y$, it is invertible with inverse norm $O(1/\log(1/h))$.
At a real shift the odd block itself is positive definite.
Its Schur complement on the even
coordinates tends to $K_s$; the source columns converge in $L^2$.

For nonreal $\lambda$, the Fourier resolvent identity, including its
derivatives at distinct centres, gives
\[
 \operatorname{Im}K_s(\lambda)
       =-(\operatorname{Im}\lambda)\Gamma_{s,\lambda}^*
                                          \Gamma_{s,\lambda}.
\]
The six columns are independent: application of $\Delta^2+\lambda$
produces independent point masses and second derivatives at distinct
centres. Thus $K_s$ is invertible. Schur inversion in the full nine-site
formula proves \eqref{eq:threeclusterres}, and Lemma~\ref{lem:pseudo}
shows that $R_s$ is the resolvent of a self-adjoint operator.

At a real shift with the stated inertia, the same Schur argument gives
three negative eigenvalues for the full congruent nine-site matrix.
For fixed $s,\beta,\vartheta$ and sufficiently small $h$, its scalar
inverse strengths have three negative centres and six positive outer
entries, by \eqref{eq:triscaling}. The two negative indices therefore
agree. Theorem~\ref{thm:inertia} proves the common bound and also
invertibility at this shift. The same matrix convergence proves the
last resolvent assertion.
\end{proof}

\begin{theorem}[Two independent Hessians with a common lower bound]
\label{thm:twohessians}
For every $\varpi\in\R$ there are nine-site scalar models and ordinary
smooth fixed-profile scalar potentials, all with a common spectral
lower bound, converging in norm resolvent sense to an operator
$H^{\rm tf}_\varpi$ whose active source space is exactly
\[
 \operatorname{span}\{(\partial_1^2-\partial_2^2)G_\lambda,
                                  2\partial_1\partial_2G_\lambda\}.
\]
In particular the two independent Hessian channels are simultaneous,
not limits of separate operator families. With
\begin{equation}\label{eq:tfres}
 W_\lambda=((\partial_1^2-\partial_2^2)G_\lambda,
                              2\partial_1\partial_2G_\lambda),\qquad
 a_\varpi(\lambda)=\varpi+\frac{1-\sqrt\lambda}{16},
\end{equation}
the limiting resolvent is
\[
 (H^{\rm tf}_\varpi+\lambda)^{-1}
 =R_0(\lambda)-a_\varpi(\lambda)^{-1}
                                      W_\lambda W_{\bar\lambda}^*.
\]
This holds at nonreal shifts and at positive real shifts for which
$a_\varpi(\lambda)\ne0$. A construction is obtained by taking any
$s_n\downarrow0$, tuning the three radial triplets by
\begin{equation}\label{eq:twotuning}
 \beta_s=-\tfrac12m_s(1),\qquad
 \vartheta=\tfrac38\varpi-\tfrac3{128},
\end{equation}
and then choosing their inner widths $h_n=o(s_n)$ sufficiently small.
The smoothing widths are chosen still smaller. The centre tuning is not
a fixed parameter: by Lemma~\ref{lem:radialcancel},
\begin{equation}\label{eq:betadiverges}
 \beta_s=\frac1{8\pi}\Bigl(\log\frac1s-\tfrac12\log3-\ell_1-\tfrac14\Bigr)
        +O(s^2)\longrightarrow+\infty .
\end{equation}
The theorem supplies existence of the scale choices $h_n$ and of the
smoothing widths, not a sharp relation between them and $s_n$.
\end{theorem}
\begin{proof}
\emph{The two symmetry blocks.}
Write $q_s=-2\beta_s+\log\lambda/(16\pi)$ and
$t=\vartheta+3(1-\sqrt\lambda)/64$.
Every three-by-three block of $K_s$ has a constant diagonal and a
constant off-diagonal. The normalized constant vector and its
two-dimensional orthogonal complement therefore reduce $K_s$ into
one copy of $B_A$ and two identical copies of $B_E$, where
\begin{align}
 B_A&=\begin{pmatrix}u_A&r_A\\r_A&v_A\end{pmatrix},
 &(u_A,r_A,v_A)&=(g+2f,\ q_s+2m_s,\ t+2n_s),\nonumber\\
 B_E&=\begin{pmatrix}u_E&r_E\\r_E&v_E\end{pmatrix},
 &(u_E,r_E,v_E)&=(g-f,\ q_s-m_s,\ t-n_s).
                                                        \label{eq:twoblocks}
\end{align}
Here $f=f_\lambda(r)$ and all cross entries have shift $\lambda$.
Lemma~\ref{lem:radialcancel} and \eqref{eq:twotuning} give, with
$L=\log(1/s)$,
\begin{align}
 u_A&\longrightarrow3g_\lambda,&
 r_A&=-6c_0L+O(1),&v_A&=O(1),\nonumber\\
 u_E&\sim3c_0s^2L,&r_E&=O(s^2),&
 v_E&\longrightarrow\tfrac38 a_\varpi(\lambda).
                                                        \label{eq:twoasym}
\end{align}
Consequently $B_A^{-1}\to0$. For nonreal $\lambda$, or at a fixed
real shift with $a_\varpi(\lambda)\ne0$, direct two-by-two inversion gives
\begin{equation}\label{eq:twoinverse}
 (B_E^{-1})_{11}=O((s^2L)^{-1}),\qquad
 (B_E^{-1})_{12}=O(L^{-1}),\qquad
 (B_E^{-1})_{22}\longrightarrow\frac8{3a_\varpi(\lambda)}.
\end{equation}

\emph{Convergence of the actual source operators.}
All constant-coordinate source columns remain bounded in $L^2$, so
their contribution tends to zero with $B_A^{-1}$. If $z\in\R^3$ has
$\sum z_j=0$, then translation differentiability of $G_\lambda$ gives
\[
 \left\|\sum_jz_jG_\lambda(\cdot-sa_j)\right\|_2=O(s)|z|.
\]
The corresponding Hessian column tends in $L^2$ to
$\sum z_j\partial_{a_j}^2G_\lambda$, by continuity of translations
on $L^2$. Thus the first entry in \eqref{eq:twoinverse} contributes
$O(1/L)$ to the resolvent; its off-diagonal contribution also vanishes.
For the orthonormal zero-sum vectors
\[
 z^{(1)}=(2,-1,-1)/\sqrt6,\qquad
 z^{(2)}=(0,1,-1)/\sqrt2,
\]
the two limiting Hessian columns are
\[
 \sqrt{3/8}(\partial_1^2-\partial_2^2)G_\lambda,
 \qquad-\sqrt{3/8}\,2\partial_1\partial_2G_\lambda.
\]
Combining this factor with the last entry of \eqref{eq:twoinverse}
proves that $R_s(\lambda)$ converges to the claimed formula, which is
a resolvent of a self-adjoint operator by Lemma~\ref{lem:pseudo},
applied to the operators of Lemma~\ref{lem:separatedtriplets}. The two columns of $W$ are
independent by their Fourier polynomials, and the scalar inverse in
\eqref{eq:tfres} is nonzero. The correction has exactly rank two with
the asserted active source space.

\emph{A single spectral lower bound, before taking either limit.}
Choose a real $\lambda_0>0$ so large that
$a_\varpi(\lambda_0)<0$. Then $u_E>0$ and $v_E<0$ for all sufficiently
small $s$, so $\det B_E<0$. Also $u_A\to3g_{\lambda_0}>0$ and
$u_Av_A-r_A^2<0$ eventually. Each real two-by-two block has one
positive and one negative eigenvalue. Hence $K_s(\lambda_0)$ has
exactly three negative and three positive eigenvalues, for all small
$s$. This fixed shift does not depend on $s$.

For each such $s_n$, Lemma~\ref{lem:separatedtriplets} permits a choice
$0<h_n<s_n/n$ for which the finite nine-site matrix has the same
three negative eigenvalues and its resolvent at $\lambda_0$ differs
from $R_{s_n}(\lambda_0)$ by less than $1/n$. Choose it also small
enough that $b_c<0<b_o$. Thus $H_{s_n,h_n}\ge-\lambda_0$ for
every $n$ after deleting finitely many initial terms. The same inner
choice can make the resolvent error at $\lambda=i$ less than $1/n$.
The already proved outer limit then proves norm resolvent convergence
of this one sequence of nine-site models. There is no interchange of
two nonuniform limits.

\emph{Ordinary local potentials.}
For each fixed member of this sequence its shifted form, at
$\lambda_0+1$, obeys
\[
 q_{s_n,h_n}[u]+(\lambda_0+1)\|u\|_2^2
       \ge c_n\|u\|_{H^2}^2,\qquad c_n>0.
\]
Indeed the fixed-configuration relative form matrix is positive
definite at this shift and equals the identity off a finite-dimensional
space. Replace each point term by
$\alpha_j\int\rho_\varepsilon(x-x_j)|u(x)|^2\,dx$, with the same
fixed profile from Section~\ref{sec:loc-smooth}.
Estimate \eqref{eq:mollify} bounds the
unsigned difference by
$C\varepsilon^{1/2}\sum_j|\alpha_j|\|u\|_{H^2}^2$.
Choose $\varepsilon_n<h_n/n$ so small that this is at most
$\eta_n$ times the shifted point form, where $\eta_n\downarrow0$
and $\eta_n<1/2$. The perturbed shifted form stays nonnegative.
Since the point shifted operator is at least $I$, the form perturbation
estimate \eqref{eq:formpert} gives resolvent error at most
$\eta_n/(1-\eta_n)$. Thus the ordinary smooth multiplication
potentials have the common lower bound $-\lambda_0-1$ and the same
norm resolvent limit. Their operator domains are $H^4$ for every $n$.
\end{proof}

This settles stable simultaneous activation, with nine scalar sites
and three attractive sites. It does not prove that nine or three is
minimal. In particular, asymmetric five-site sufficiency and all
intermediate site counts remain open. The construction is genuinely
multiscale and is not a counterexample to the fixed-geometry lower
bound in Theorem~\ref{thm:fourlower}. Nor does it supply arbitrary
six-channel coupling matrices: here only the two traceless Hessians
are active and they have one common real parameter.

\subsection{Geometry-independent bounds and an optimal dipole construction}\label{sec:loc-minimum}

\label{sec:generalminimality}
We now allow arbitrary distinct sites $x_j(h)\to0$, including
degenerating configurations and unrelated collision scales. The number
of sites is bounded by a fixed integer $N$. At a fixed nonreal shift
write $\mathcal T_h=R_h(\lambda)-R_0(\lambda)$ and let $\mathcal T$ be its norm limit.
By Lemma~\ref{lem:sixjets}, the range of $\mathcal T$ consists of fields
\[
 \widehat u(\xi)=\frac{P(\xi)}{|\xi|^4+\lambda},\qquad \deg P\le2.
\]
Let $\mathcal P_{\mathcal T}$ be this polynomial space and define its derivative
closure by
\[
 \mathcal D(\mathcal P_{\mathcal T})
 =\operatorname{span}\{\partial^\alpha P:P\in\mathcal P_{\mathcal T},
                                      \ |\alpha|\le2\}.
\]
The dimension of this closure can exceed the rank of the resolvent
correction. It gives a geometric bound on the microscopic site count.
For a fixed finite node set, de Boor and Ron form an interpolation
space from the lowest-degree
homogeneous terms of exponential sums and show that it is closed
under differentiation \cite[Section~3]{dBR90}. We use the moving
site exponentials to bound the dimension for a norm limit of resolvents.

\begin{theorem}[A site bound valid at arbitrary collision scales]
\label{thm:derivativeclosure}
Every norm resolvent limit of at most $N$ scalar point sites converging
to zero satisfies
\begin{equation}\label{eq:derivativeclosure}
                  \dim\mathcal D(\mathcal P_{\mathcal T})\le N.
\end{equation}
No common lower bound or fixed reference geometry is required.
In particular, two independent Hessian directions require at least
five sites, and three independent Hessian directions require at least
six. A channel with one nonzero quadratic symbol $q$ requires at least
$2+\rank(\nabla^2q)$ sites. Two independent dipole directions require
at least three sites. The conclusions also apply to smooth families
whose resolvents are norm close to the corresponding point models.
\end{theorem}
\begin{proof}
Put $\mathcal S_h=\operatorname{span}\{G_\lambda(\cdot-x_j(h))\}$,
so $\dim\mathcal S_h\le N$ and $\operatorname{ran}\mathcal T_h\subset\mathcal S_h$.
Choose $\chi\in C_c^\infty(\R^2)$ equal to one on a neighbourhood
containing every site for all sufficiently small $h$. The map
$\Delta^2+\lambda:L^2\to H^{-4}$ is an isomorphism: its Fourier
symbol is comparable in modulus to $(1+|\xi|^2)^2$.
Multiplication by $x_i\chi$ is bounded on $H^{-4}$, by duality and
the product rule on $H^4$. Hence
\[
 \mathcal A_i=R_0(\lambda)(x_i\chi)(\Delta^2+\lambda),\qquad i=1,2,
\]
is bounded on $L^2$. It preserves $\mathcal S_h$, since it multiplies
the charge at $x_j(h)$ by its $i$th coordinate.

For any $u=\mathcal Tf$, take $u_h=\mathcal T_hf\to u$. Every finite word in
$\mathcal A_1,\mathcal A_2$ applied to $u_h$ stays in $\mathcal S_h$
and converges to that word applied to $u$. At the limiting point
source, multiplication by $x_i$ differentiates its Fourier polynomial
by $i\partial_{\xi_i}$; the cutoff contributes nothing because it is
identically one near zero. Thus all fields belonging to
$\mathcal D(\mathcal P_{\mathcal T})$ have convergent approximants in
$\mathcal S_h$. If more than $N$ were independent, the Gram matrices
of a finite independent set and of its approximants would converge to
a positive definite matrix, contradicting $\dim\mathcal S_h\le N$.
This proves \eqref{eq:derivativeclosure}.

For two independent homogeneous quadratic parts, their first
derivatives span both linear variables. Otherwise both symmetric
Hessian matrices would have range in one common line, making both
quadratics multiples of the square of the same linear form. Their
second derivatives also supply the constant. Counting degrees gives
$2+2+1=5$. With three quadratic parts the count is $3+2+1=6$.
For a single quadratic the first-derivative span has dimension
$\rank(\nabla^2q)$, and the constant supplies one further dimension.
The dipole assertion similarly counts two linear terms and a constant.
Norm close smooth families have exactly the same limiting correction,
so inherit the conclusions.
\end{proof}

Theorem~\ref{thm:fourlower} is the fixed-configuration case. The
nine-site upper
construction still leaves a gap; this theorem is a necessary bound,
not a five-site realization.

\begin{theorem}[Six sites cannot activate all six channels]
\label{thm:sixsiteobstruction}
Let arbitrary distinct scalar sites converge to zero. If the limiting
resolvent correction has rank six, then at least seven sites are
necessary for every sufficiently small collision parameter. In
particular no family with at most six sites, with any geometry or
strengths, can converge to a full six-channel active relation.
This conclusion does not require a common spectral lower bound.
It also holds for smooth realizations that are norm resolvent close to
point models with at most six sites.
\end{theorem}
\begin{proof}
Fewer than six sites are excluded by the rank of the correction.
Suppose, along a sequence, that there are exactly six and $\rank\mathcal T=6$.
Use the six multiindices $|\alpha|\le2$, put
$\gamma_\alpha=(-1)^{|\alpha|}\partial^\alpha G_\lambda$, and let
$\Gamma$ have these six columns. Take test functions
$\phi_\alpha(x)=\chi(x)x^\alpha/\alpha!$ with a common cutoff
equal to one near zero. The complex linear map
\[
 Fu=\big(\langle(\Delta^2+\lambda)u,\phi_\alpha\rangle\big)_\alpha
\]
is bounded on $L^2$ and satisfies $F\Gamma=I$. Here the pairing is
distributional, without complex conjugation. On the six translates,
\[
 C_h:=F\Phi_{h,\lambda}
       =\big(x_j(h)^\alpha/\alpha!\big)_{\alpha,j}.
\]

Choose $f_1,\ldots,f_6$ so that $\mathcal Tf_j$ form a basis of
$\operatorname{ran}\Gamma$, and put $X_h=(\mathcal T_hf_1,\ldots,\mathcal T_hf_6)$.
Then $X_h\to X$, $FX$ is invertible, and the six columns of $X_h$
span $\mathcal S_h$ for small $h$. Consequently $F|_{\mathcal S_h}$
and $C_h$ are invertible, and
\[
 W_h:=\Phi_{h,\lambda}C_h^{-1}
     =X_h(FX_h)^{-1}\longrightarrow\Gamma.
\]
This argument obtains the required quadratic interpolation basis even
for a degenerating geometry; it does not assume an interpolation
condition in advance. With $Q_h=B_h+M_{h,\lambda}$, set
\[
 A_h=C_h^{-T}Q_hC_h^{-1},\qquad
 \mathcal T_h=-W_hA_h^{-1}W_{h,\bar\lambda}^*.
\]
The same construction at $\bar\lambda$ gives convergence of the
second source map. Continuous left inverses recover a limit of
$A_h^{-1}$. Its rank is six because $\rank\mathcal T=6$; hence $A_h$
converges to an invertible matrix and is uniformly bounded.

Let $c(x)=(x^\alpha/\alpha!)_\alpha$. The exact equality
$Q_h=C_h^TA_hC_h$ gives, whenever $i\ne j$,
\[
 G_\lambda(x_i-x_j)=c(x_i)^TA_hc(x_j).
\]
The right side is a polynomial in two variables with uniformly bounded
mixed second derivatives near zero. For four distinct indices and
$a=x_i-x_k$, $b=x_j-x_l$, it follows that
\begin{equation}\label{eq:rectupper}
 |G_\lambda(x_i-x_j)-G_\lambda(x_i-x_l)
       -G_\lambda(x_k-x_j)+G_\lambda(x_k-x_l)|\le C|a||b|.
\end{equation}

Among any six distinct planar points there are two disjoint segments
with $|a\cdot b|\ge|a||b|/2$. To see this, partition unoriented
directions into three intervals of width $\pi/3$ and color the edges
of the complete six-vertex graph accordingly. Suppose that no color
contains two disjoint edges. In a given color choose two edges $uv,uw$
if such a pair exists. Every edge of that color either meets $u$ or is
$vw$. If $vw$ occurs, all edges of that color lie in the triangle
$uvw$. Otherwise they share $u$. Each color therefore has at most five
edges, with equality only for a full star. Covering all fifteen edges
would require three full stars, but any two full stars share an edge. Orient one of the
segments so that $a\cdot b\ge|a||b|/2$.

Write $\delta_h=\max_{r,t}|x_r-x_t|\to0$. The local expansion
\eqref{eq:local} gives
\[
 \nabla^2G_\lambda(z)=2c_0\log|z|\,I+E_\lambda(z),
 \qquad \|E_\lambda(z)\|\le C_\lambda,
 \qquad c_0=1/(8\pi).
\]
Twice using the fundamental theorem of calculus represents the
rectangular difference in \eqref{eq:rectupper} as
\[
 -\int_0^1\!\int_0^1
 a^T\nabla^2G_\lambda(x_k+ta-x_l-ub)b\,dt\,du.
\]
This identity also holds for intersecting or overlapping segments.
The kernel is $C^1$ at zero. On a line through zero its gradient
is absolutely continuous on bounded intervals, since its derivative
is $O(1+|\log|z||)$. For fixed nonzero $a,b$, the logarithm is
integrable on the parameter square: an invertible affine change of
variables gives a locally integrable planar logarithm when $a,b$
are independent. When they are parallel, integration along a
parameter with nonzero slope reduces the singularity to the
one-dimensional integral of $|\log|t||$. Thus the two integrations
are justified separately for each configuration. Every distance in this integral is at most $\delta_h$.
Taking real parts therefore bounds that difference below by
\[
 c_0|a||b|\log(1/\delta_h)-C_\lambda|a||b|,
\]
contradicting \eqref{eq:rectupper}. The norm close smooth case again
has the same limiting correction.
\end{proof}

\begin{theorem}[Both dipoles from an obtuse three-site configuration]
\label{thm:obtusedipoles}
For every real symmetric $2\times2$ matrix $S$, there are scalar
models with exactly three sites and two attractive strengths,
and ordinary smooth fixed-profile realizations with a common lower
bound, whose norm resolvent limit has active source space
\[
 \operatorname{span}\{G_\lambda,D_{1,\lambda},D_{2,\lambda}\},
 \qquad D_{i,\lambda}=-\partial_iG_\lambda.
\]
Its resolvent is
\begin{equation}\label{eq:obtuseres}
 R_S(\lambda)=R_0(\lambda)-g_\lambda^{-1}G_\lambda G_{\bar\lambda}^*
   -D_\lambda\left(S-\frac{\log\lambda}{16\pi}I\right)^{-1}
                                                    D_{\bar\lambda}^*,
 \qquad D_\lambda=(D_{1,\lambda},D_{2,\lambda}).
\end{equation}
This formula holds at nonreal shifts and at positive real $\lambda$
with $\log\lambda/(16\pi)\notin\sigma(S)$. The monopole is pinned: in
the chart
$K_\lambda=\Theta+\diag\bigl(g_\lambda,-\tfrac{\log\lambda}{16\pi},
-\tfrac{\log\lambda}{16\pi}\bigr)$, which extends \eqref{eq:Theta}, the
coupling matrix is $\Theta=\diag(0,S)$. Thus the monopole--dipole
couplings vanish, and every function in the range of $R_S(\lambda)$
vanishes at the origin. Three sites and two
attractive sites are minimal for simultaneous activation of both
dipoles in the point model and its norm resolvent close smooth realizations. This is a three-channel family with
arbitrary real dipole block $S$; it is not a classification of all
monopole--dipole couplings or of the six-channel relations.
\end{theorem}
\begin{proof}
Take the obtuse triangle
\[
 x_0=0,\qquad x_1=(h,0),\qquad x_2=(h/2,h/4),\qquad
 V=\begin{pmatrix}-1&-2\\1&-2\\0&4\end{pmatrix}.
\]
Its affine moment matrix satisfies
$\boldsymbol1^TV=0$ and $(x_0,x_1,x_2)V=hI$.
Let $M_{h,1}=(G_1(x_i-x_j))_{i,j=0}^2$ and put
\begin{align}
 Y_h&=h^2S-V^TM_{h,1}V,\nonumber\\
 b_0&=\tfrac12Y_{11}+\tfrac14Y_{12},\qquad
 b_1=\tfrac12Y_{11}-\tfrac14Y_{12},\qquad
 b_2=\tfrac1{16}(Y_{22}-4Y_{11}).\label{eq:obtusetuning}
\end{align}
These are exact inverse strengths. Indeed direct multiplication gives
\[
 V^T\diag(b_0,b_1,b_2)V
 =\begin{pmatrix}b_0+b_1&2(b_0-b_1)\\
                 2(b_0-b_1)&4(b_0+b_1)+16b_2\end{pmatrix}=Y_h.
\]
By \eqref{eq:local} and the affine moments,
$V^TM_{h,1}V=2c_0h^2\log(1/h)I+O(h^2)$, so
\begin{equation}\label{eq:obtusesigns}
 \frac{(b_0,b_1,b_2)}{c_0h^2\log(1/h)}
                   \longrightarrow(-1,-1,3/8).
\end{equation}
There are thus exactly two attractive sites for every sufficiently
small $h$.

Use the source change
\[
 C_h=\begin{pmatrix}1&1&1\\0&h&h/2\\0&0&h/4\end{pmatrix},
 \qquad W_{h,\lambda}=\Phi_{h,\lambda}C_h^{-1}.
\]
The columns of $C_h^{-1}$ are $e_0$ and the two columns of $V/h$.
Consequently $W_{h,\lambda}\to(G_\lambda,D_{1,\lambda},D_{2,\lambda})$
in operator norm, by the $L^2$ translation Taylor estimate.
In the congruent matrix $K_{h,\lambda}=C_h^{-T}Q_hC_h^{-1}$,
the monopole entry tends to $g_\lambda$, the two mixed entries are
$O(h\log(1/h))$, and the dipole block is exactly
\[
 S+h^{-2}V^T(M_{h,\lambda}-M_{h,1})V
       \longrightarrow S-\frac{\log\lambda}{16\pi}I.
\]
The last limit is \eqref{eq:I2} in the two coordinate directions;
the off-diagonal integral vanishes by symmetry. Equivalently it follows
by Taylor expanding $G_\lambda-G_1$ through degree two.
For nonreal $\lambda$ the limiting blocks are invertible. Substitution
in \eqref{eq:res} proves \eqref{eq:obtuseres}, a self-adjoint resolvent
by Lemma~\ref{lem:pseudo}. The correction has rank three and the
displayed active space. Each $D_{i,\lambda}$ is continuous and vanishes
at the origin, by termwise differentiation of \eqref{eq:local}, and
$G_{\bar\lambda}^*f=(R_0(\lambda)f)(0)$. Evaluating \eqref{eq:obtuseres}
at the origin therefore gives $(R_S(\lambda)f)(0)=0$ for every $f$.

Choose a real $\lambda_0$ with
$\log\lambda_0/(16\pi)>\lambda_{\max}(S)$.
The limiting matrix has one positive and two negative eigenvalues.
For all small $h$ its congruent $Q_h(\lambda_0)$ therefore has
exactly two negative eigenvalues and no zero eigenvalue. This matches
$n_-(B_h)=2$ from \eqref{eq:obtusesigns}. Theorem~\ref{thm:inertia}
gives the common lower bound $H_h\ge-\lambda_0$.
At the larger shift $\lambda_0+1$ each fixed point form controls
$H^2$ with some positive constant depending on $h$. Estimate
\eqref{eq:mollify} permits widths $\varepsilon(h)=o(h)$, small enough
that the absolute smoothing error is $o(1)$ relative to this shifted
form. The form estimate \eqref{eq:formpert} then gives norm resolvent
error $o(1)$ and the common lower bound $-\lambda_0-1$ for the
ordinary smooth potentials, with the same fixed profile at all sites.
Finally Theorems~\ref{thm:derivativeclosure} and
\ref{thm:attractiverank} give the matching lower bounds of three
sites and two attractive sites.
\end{proof}

These results strengthen the necessary conditions and give a sharp
simultaneous-dipole construction. They do not prove five-site
simultaneous Hessians, the optimal Hessian site or attractive-site
counts, or sufficiency of seven or more sites for every real
six-channel relation. Those are separate unresolved construction and
classification questions.

\appendix

\section{Admissible jets and the missing-multipole obstruction}\label{fd:sec:found-restrictions}

We use $\mathcal R_0(z)=(A_0-z)^{-1}$ and denote the graph trace,
defect map and Weyl matrix by $\tau$, $\gamma_z$ and $M_z$, respectively.

\subsection{The point restriction and its defect space}\label{fd:sec:found-point}

Put $P=(-\Delta)^m$, $\Dom(P)=H^{2m}(\R^n)$, and let $S$ be the closure
of its restriction to $C_c^\infty(\R^n\setminus\{0\})$.
Use the Fourier convention of Section~\ref{sec:loc-model} in dimension $n$.
Write
\begin{equation}\label{fd:eq:count}
t=2m-\frac n2,\qquad
 k=\max\{-1,\lceil t\rceil-1\},\qquad
 r=\binom{n+k}{n},
\end{equation}
with $\binom{n-1}{n}=0$ and an empty jet list when $k=-1$.

\begin{lemma}[A point-supported distribution and its Sobolev order]
\label{fd:lem:point}
A distribution supported at $0$ is a finite sum of derivatives of
$\delta_0$. A nonzero sum with largest derivative degree $d$ belongs to
$H^{-s}(\R^n)$ precisely when $d<s-n/2$.
\end{lemma}
\begin{proof}
The structure theorem is \cite[Theorem~2.3.4]{Hor}. The Fourier
transform is a polynomial with nonzero leading homogeneous part. Its
squared $H^{-s}$ norm at infinity is comparable to
$\int_1^\infty r^{2d-2s+n-1}\,\dd r$, giving the strict threshold.
\end{proof}

The following Sobolev domain description gives the admissible multipoles
of Karpeshina and Pavlov \cite[Sections~1--2]{KP}.
\begin{theorem}[Admissible jets and deficiency indices]\label{fd:thm:defect}
The closed restriction has domain
\begin{equation}\label{fd:eq:old2}
\Dom(S)=\{u\in H^{2m}(\R^n):\partial^\alpha u(0)=0,
                                     \ |\alpha|\le k\}.
\end{equation}
For $z\notin\R$, its defect space is
\begin{equation}\label{fd:eq:old3}
\ker(S^*-z)=\operatorname{span}\left\{
 \mathcal F^{-1}\frac{(i\xi)^\alpha}{|\xi|^{2m}-z}:
                                      |\alpha|\le k\right\}.
\end{equation}
Its deficiency indices are $(r,r)$. If $n\ge4m$, the point restriction
is essentially self-adjoint.
\end{theorem}
\begin{proof}
In the graph norm $H^{2m}$ the annihilator of the punctured test
functions consists of the $H^{-2m}$ distributions supported at zero.
Lemma~\ref{fd:lem:point} identifies it with the jets in \eqref{fd:eq:old2}.
The defect equation makes $(P-z)u$ a distribution of the same type.
Dividing its Fourier polynomial by $|\xi|^{2m}-z$ gives
\eqref{fd:eq:old3}, with one independent field per admissible monomial
\cite[Sections~1--2]{KP}.
\end{proof}

\begin{corollary}[Several distinct points]
For $N$ fixed distinct punctures the index is $(Nr,Nr)$, and the minimal
domain kills the admissible jets at every puncture.
\end{corollary}
Disjoint cutoffs reduce this statement to the one-point result. For
$m=2$ the indices per point in dimensions $1,2,3,4$ are respectively
$4,6,10,5$. The planar channels are listed in Lemma~\ref{lem:sixjets}.

\subsection{A missing-multipole obstruction}\label{fd:sec:found-obstruction}

\label{fd:sec:counterexample}
We examine uniqueness in the ambient space defined in \cite{KT}.
In this section
$\Omega=B(0,1)\subset\R^2$, $m=2$, and the outer conditions are the
Navier conditions $u=\Delta u=0$ on $r=1$.
For $i=0,1$ the six functionals printed in
\cite[Section~2, pp.~340--341]{KT} have the form
\[
 \mathfrak f_i^0(u)=-\lim_{r\downarrow0}\int_{|x|=r}
                   \partial_\nu\Delta^{1-i}u\,\dd S,
 \quad
 \mathfrak f_i^s(u)=\lim_{r\downarrow0}\int_{|x|=r}
                   \omega_s\Delta^{1-i}u\,\dd S,
 \quad s=1,2.
\]
The following zero identities hold for either choice of radial normal.
The singular ambient space on those pages requires these limits to be
finite and the punctured expression
$\Delta^2u-\sum_j\mathfrak f_j(u)\widetilde G_j$ to belong to $L^2$.
It does not impose a bound excluding higher angular multipoles.

\begin{proposition}[An explicit $L^2$ counterexample to the printed uniqueness]
\label{fd:prop:counterexample}
The function
\begin{equation}\label{fd:eq:counterexample}
 u(r,\theta)=(r^4-4r^2+3)\cos(2\theta),\qquad 0<r<1,
\end{equation}
belongs to that singular ambient space and to $L^2(\Omega)$, is nonzero,
and satisfies
\[
 \Delta^2u=0\quad\hbox{on }\Omega\setminus\{0\},\qquad
 u|_{r=1}=\Delta u|_{r=1}=0,\qquad
 \mathfrak f_i^s(u)=0\quad(i=0,1;\ s=0,1,2).
\]
Thus the uniqueness assertion of \cite[Theorem~2.1]{KT} is
false in its printed ambient space for $m=n=2$.
The same example satisfies the homogeneous conditions of its
Theorem~2.2 for every set of bounded linear functionals $c_j$.
Restricting that ambient space only by $u\in L^2$ does not repair the
uniqueness assertion.
More precisely, let $\mathcal K$ consist of $L^2(B)$ fields which are
smooth off zero, are $H^4$ in an outer collar, solve the homogeneous
punctured equation, satisfy the Navier traces and have all six printed
fluxes zero. Then
\[
 \mathcal K=\operatorname{span}\{p(r)\cos2\theta,p(r)\sin2\theta\},
 \qquad p(r)=r^4-4r^2+3.
\]
\end{proposition}
\begin{proof}
For every real exponent $a$ and integer $\ell$,
\[
 \Delta(r^a\cos\ell\theta)=(a^2-\ell^2)r^{a-2}\cos\ell\theta.
\]
Hence $\Delta u=12(r^2-r^{-2})\cos2\theta$ and $\Delta^2u=0$
off zero. Both outer traces vanish. The function is smooth on every
compact subset of the punctured disk and is bounded. In fact
\[
 \norm u_2^2=\pi\int_0^1(r^4-4r^2+3)^2r\,\dd r=\frac{19\pi}{15}>0.
\]
Each integral defining a zeroth flux has angular factor $\cos2\theta$;
each first moment has factor $\cos\theta\cos2\theta$ or
$\sin\theta\cos2\theta$. All three angular integrals are zero at every
positive radius. Thus all six limits exist and equal zero. The remaining
ambient-space expression is the zero function, which is in $L^2$.
For any linear $c_j$, $c_j(\Delta^2u)=c_j(0)=0$, proving the assertion
about Theorem~2.2.

The missing source is visible distributionally. With
$\Phi=r^2\log r/(8\pi)$ one has
$(\partial_1^2-\partial_2^2)\Phi=\cos2\theta/(4\pi)$, and therefore
\begin{equation}\label{fd:eq:counterexamplesource}
 \Delta^2u=12\pi(\partial_1^2-\partial_2^2)\delta_0
 \quad\hbox{in }\mathcal D'(\Omega).
\end{equation}
The other two terms of \eqref{fd:eq:counterexample} are biharmonic
polynomials. The listed radial and first angular moments miss this
trace-free Hessian source.

To prove the last assertion, extend any $v\in\mathcal K$ as its
$L^2$ distribution across zero. The distribution $\Delta^2v$ is
supported there. After localization it belongs to $H^{-4}$, so
Lemma~\ref{fd:lem:point} bounds its derivative degree by two. The six
zero-energy Navier Green fields with sources
$\delta_0,-\partial_1\delta_0,-\partial_2\delta_0,
\Delta\delta_0,(\partial_1^2-\partial_2^2)\delta_0,
\partial_1\partial_2\delta_0$
therefore span all these homogeneous $L^2$ Navier solutions.
Indeed, subtract the Green field with the same source. The difference
is biharmonic in the distributional sense across zero and hence smooth
there by interior elliptic regularity \cite[Theorem~6.29]{Grubb}. It is $H^4$
in the collar and has zero Navier traces, hence is in the kernel of
the invertible Navier realization and vanishes.
The local flux integrations in Proposition~\ref{fd:prop:flux}, proved
independently below, distinguish the first four sources: the radial
fluxes recover the monopole and trace-Hessian coefficients, and the
unscaled first moments of the two source derivatives
$\partial_{y_t}G$, whose sources are $-\partial_t\delta_0$, are
$-\delta_{st}/2$. Their remaining printed fluxes are zero. The last two Green
fields have angular order two on the centered disk and all the
printed fluxes vanish. Thus precisely their coefficients are free.
Formula~\eqref{fd:eq:counterexamplesource} identifies the first field;
for the second,
\[
 \Delta^2\bigl(p(r)\sin2\theta\bigr)=24\pi\partial_1\partial_2\delta_0.
\]
This proves the exact two-dimensional assertion. The collar assumption
makes the outer traces and the use of the fixed Navier realization
explicit; the original counterexample is smooth in that collar.
\end{proof}

\begin{proposition}[Infinite ambiguity without an $L^2$ growth condition]
\label{fd:prop:invisibleinfinite}
In the same printed ambient space the homogeneous Navier problem has
an infinite-dimensional kernel on which all six functionals vanish.
\end{proposition}
\begin{proof}
For every integer $\ell\ge2$ put
$v_\ell=(r^\ell-r^{-\ell})\cos\ell\theta$.
These functions are harmonic on the punctured disk, vanish on its outer
circle, and have zero zeroth and first angular moments. Thus they meet
every defining condition of the ambient space, with all fluxes and the
punctured $\Delta^2v_\ell$ equal to zero. Their distinct angular modes
are linearly independent. They are not $L^2$; the preceding proposition
supplies the separate $L^2$ obstruction.
\end{proof}

The specific missing step is the assertion in
\cite[p.~347]{KT} that finiteness of the displayed zeroth and first
moments makes every positive Taylor-weighted moment vanish. For
\eqref{fd:eq:counterexample} the outward-ball normal gives instead
\[
 \int_{|x|=r}\partial_r\Delta u\,x_1^2\,\dd S
       =12\pi(1+r^4)\longrightarrow12\pi.
\]
Higher angular data therefore cannot be dropped from that Green formula.
If its ambient-space differential expression were interpreted
distributionally on the whole disk, the point distributions in
\eqref{fd:eq:counterexamplesource} would appear. But that interpretation
would also exclude the advertised Green channels: $\Delta^2G=\delta_0$,
and subtracting ordinary Green functions does not cancel $\delta_0$.
It does not supply a consistent alternative reading of the printed
definition. A repaired direct-sum domain is constructed below.

\Needspace{8\baselineskip}

\section{Singular coordinates and flux normalization}\label{fd:sec:found-coordinates}

\subsection{Fundamental solutions and flux coordinates}\label{fd:sec:found-flux}

Let $\sigma_{n-1}=|\mathbb S^{n-1}|$.

\begin{proposition}[The exceptional logarithmic dimensions]
\label{fd:prop:fundamental}
In the normalization of \cite[Section~2.6, formula~(2.62)]{GGS},
for even $n\ge2$ and $k=m-n/2\ge0$, a fundamental solution of
$(-\Delta)^m$ is
\begin{equation}\label{fd:eq:old4}
\Phi_{m,n}(x)=
 \frac{(-1)^{k+1}|x|^{2k}\log|x|}
 {2^{2m-1}\pi^{n/2}(m-1)!\,k!}.
\end{equation}
Away from these exceptional dimensions, including $n=1$, it is
\begin{equation}\label{fd:eq:old5}
\Phi_{m,n}(x)=\frac{\Gamma(n/2-m)}{4^m\pi^{n/2}(m-1)!}
                                     |x|^{2m-n}.
\end{equation}
For $\Delta^m$ multiply these expressions by $(-1)^m$. In particular
$\Phi_{2,2}=|x|^2\log|x|/(8\pi)$ with either convention.
\end{proposition}

For the next proposition use the $\Delta^m$ fundamental solution
$\Psi=(-1)^m\Phi_{m,n}$ and a Green kernel
$G(x,y)=\Psi(x-y)+H(x,y)$, with $H$ smooth jointly near $(0,0)$.
Define $U_j=\Delta_y^jG(\cdot,0)$ and
$V_{j,s}=\partial_{y_s}\Delta_y^jG(\cdot,0)$, $0\le j<m$.

\begin{proposition}[Corrected local flux duality]\label{fd:prop:flux}
For these channels plus functions smooth near zero, put
\begin{align*}
 F_i^0(u)&=\lim_{r\downarrow0}\int_{|x|=r}
                       \partial_r\Delta^{m-i-1}u\dd S,\\
 F_i^s(u)&=-n\lim_{r\downarrow0}\int_{|x|=r}
                       \frac{x_s}{r}\Delta^{m-i-1}u\dd S.
\end{align*}
Then
\begin{equation}\label{fd:eq:old7}
F_i^0(U_j)=\delta_{ij},\quad F_i^s(V_{j,t})=\delta_{ij}\delta_{st},
 \quad F_i^s(U_j)=F_i^0(V_{j,t})=0.
\end{equation}
These identities remain valid in the logarithmic cases.
\end{proposition}
\begin{proof}
Smooth summands contribute zero in each limit. For the first functional
their constant vector flux integrates to zero and their remainder is
$O(r^n)$; for the second, the constant spherical moment is zero and
the remainder is $O(r^n)$.
For the singular summands $\Delta_y$ can be replaced by $\Delta_x$.
If $j>i$, at least $m$ Laplacians have been applied and the result is
zero off the pole. If $j=i$, $\Delta^{m-1}\Psi$ is the Laplace
fundamental solution up to an irrelevant smooth term, with gradient
$x/(\sigma_{n-1}r^n)$. Its radial flux is one.
Its source derivative is the negative spatial derivative, so
\[
 \int_{|x|=r}\frac{x_s}{r}
       \partial_{y_t}\Delta_x^{m-1}\Psi(x-y)|_{y=0}\dd S
       =-\frac1{\sigma_{n-1}}\int_{\mathbb S^{n-1}}
                                          \omega_s\omega_t\dd\omega
       =-\frac{\delta_{st}}n.
\]
Reflection symmetry kills the off-diagonal moments, and summing the
diagonal moments gives one. This explains the required factor $-n$.
If $j<i$, the corresponding radial fundamental solution has order
$2(i+1-j)-n$, possibly with a logarithm. Its radial flux and the
source-derivative first moment are $O(r^{2(i-j)}(1+|\log r|))$,
and tend to zero. The cross identities vanish identically for singular
terms by oddness. This proves all four identities without replacing a
logarithmic fundamental solution by a pure power.
\end{proof}

\paragraph{Normals and the printed normalization.}
The local singular calculations in \cite[Lemmas~3.1--3.2]{KT} use a
pure-power expression. The exceptional dimensions require the logarithmic
kernel \eqref{fd:eq:old4}. Our radial normal points out
of the small ball. With the normal pointing out of the punctured domain,
$\nu=-\partial_r$, the zeroth functional $-\int\partial_\nu$ equals
$F_i^0$. With the outward-ball normal it equals $-F_i^0$.
The latter radial convention is explicitly printed in \cite[p.~346]{KT};
it must be reconciled with the sign used in its duality claim.
For the source derivative convention used here the unscaled first moment
in that paper is always $-F_i^s/n$, not $F_i^s$. In particular its full
Kronecker duality statement requires this normalization correction even
when logarithms are absent. None of these sign choices changes the zero
flux counterexample of Proposition~\ref{fd:prop:counterexample}.
The proof above needs no equality of $\Delta_xG$ and $\Delta_yG$ for
the smooth boundary correction.

\subsection{Radial--gradient and complete jet spaces}\label{fd:sec:found-channels}

\begin{proposition}[Exact admissibility test for the selected channels]
\label{fd:prop:channels}
$U_j$ is locally $L^2$ precisely when $2j<2m-n/2$;
$V_{j,s}$ is locally $L^2$ precisely when $2j+1<2m-n/2$.
For $m=2,n=2$ this retains
$G,\partial_{y_1}G,\partial_{y_2}G,\Delta_yG$ and excludes
$\partial_{y_s}\Delta_yG$.
The full six-dimensional space has two additional trace-free Hessian
directions, for example $\partial_1^2G-\partial_2^2G$ and
$\partial_1\partial_2G$.
\end{proposition}
\begin{proof}
Locally $Pu=T+f$, where $f$ is smooth and the source $T$ has degree
$d=2j$ or $2j+1$. If $u$ is locally $L^2$, localization gives
$T\in H^{-2m}$. Conversely, for $T\in H^{-2m}$ put
$w=(P-z)^{-1}T\in L^2$ at a nonreal $z$. Since $P(u-w)=f-zw$
is locally $L^2$, interior elliptic regularity
\cite[Theorem~6.29]{Grubb} gives $u-w\in H^{2m}_{\rm loc}$.
The strict threshold follows from Lemma~\ref{fd:lem:point}. In the
plane the three quadratic monomials span one trace and two trace-free
directions.
\end{proof}

For general $m,n$, the dimension retained from the selected family is
\begin{equation}\label{fd:eq:old8}
\#\{0\le j<m:2j<t\}
 +n\,\#\{0\le j<m:2j+1<t\}.
\end{equation}
It need not be $m(n+1)$ or the full index $r$. In dimension one the
selected orders are exactly $0,1,\ldots,2m-1$ and do give the full
space. The equality $6=2(2+1)$ in the planar biharmonic case is an
equality of counts before the $L^2$ test, not an identification of bases.

A reference $L^2$ resolvent acting on an excluded
$\partial_s\Delta G$ requires an extension to a larger input space.

There is also an internal domain inconsistency in
\cite[Remark~2.1]{KT}, which places every selected Green field in the
printed ambient space. For $m=n=2$ and $V=V_{1,s}$, its only nonzero
printed flux is the matching first moment, equal to $-1/2$.
Since $\Delta^2V=0$ off the pole, the expression defining that space is
$V/2$, which is not in $L^2$. Thus that membership assertion fails
for these two fields. The finite singular direct sum in
Theorem~\ref{fd:thm:repair} is a replacement domain; when it includes
these fields it is not a subspace of the printed ambient space.

\section{Finite-trace extensions with fixed outer conditions}\label{fd:sec:found-extensions}

\subsection{The finite-trace construction}\label{fd:sec:found-traces}

Let $A_0$ be self-adjoint on a Hilbert space $\mathcal H$, and let
$\tau:\Dom(A_0)\to\C^r$ be graph-continuous and onto, with
$\ker\tau$ dense in $\mathcal H$. Take real $z_*\in\rho(A_0)$ and put
\[
 R_*=(A_0-z_*)^{-1},\quad \gamma_*=(\tau R_*)^*,\quad
 S=A_0|_{\ker\tau}.
\]
For a puncture the coordinates of $\tau$ are the admissible regular
jets. In the whole-space case $\gamma_*$ is their source-derivative
Green family, with the Fourier convention above.

\begin{theorem}[Finite-trace extensions]\label{fd:thm:trace}
\label{fd:thm:projection}\label{fd:prop:coordinates}
The operator $S$ is closed, symmetric and densely defined, with
\begin{equation}\label{fd:eq:old9}
\Dom(S^*)=\Dom(A_0)\dotplus\ran\gamma_*,\qquad
 S^*(u_0+\gamma_*c)=A_0u_0+z_*\gamma_*c.
\end{equation}
Set $\Gamma_0(u_0+\gamma_*c)=c$, $\Gamma_1(u_0+\gamma_*c)=\tau u_0$ and,
for $z\in\rho(A_0)$,
\[
 \mathcal R_0(z)=(A_0-z)^{-1},\qquad
 \gamma_z=[I+(z-z_*)\mathcal R_0(z)]\gamma_*,\qquad
 M_z=\tau(\gamma_z-\gamma_*).
\]
\begin{enumerate}[label=\textup{(\roman*)}]
\item The self-adjoint extensions of $S$ correspond one-to-one to the
pairs $(\Pi,\Theta)$ of an orthogonal projection $\Pi$ in $\C^r$ and a
Hermitian operator $\Theta$ on $\ran\Pi$, through
\begin{equation}\label{fd:eq:projectiondomain}
 \Gamma_0u\in\ran\Pi,\qquad
       \Pi\Gamma_1u=\Theta\Gamma_0u,
\end{equation}
and the resolvent of the extension is
\begin{equation}\label{fd:eq:projectionresolvent}
 R_{\Pi,\Theta}(z)=\mathcal R_0(z)+\gamma_z\Pi
  (\Theta-\Pi M_z\Pi)^{-1}_{\ran\Pi}\Pi\gamma_{\bar z}^*,
       \qquad z\notin\R.
\end{equation}
The choice $\Pi=0$ is $A_0$.
\item One has
\begin{equation}\label{fd:eq:weylidentity}
 M_z-M_w=(z-w)\gamma_{\bar w}^*\gamma_z,
 \quad M_z^*=M_{\bar z},\quad
 \operatorname{Im}M_z=(\operatorname{Im}z)\gamma_z^*\gamma_z,
\end{equation}
and the Gram matrix $\gamma_z^*\gamma_z$ is positive definite.
\item If $J:\C^s\to\C^r$ is injective and $\Theta_J$ is a Hermitian
$s\times s$ matrix, then
\begin{equation}\label{fd:eq:activerelation}
 \Gamma_0u=Jc,\qquad J^*\Gamma_1u=\Theta_Jc
\end{equation}
defines a self-adjoint extension with resolvent correction
$\gamma_zJ(\Theta_J-J^*M_zJ)^{-1}J^*\gamma_{\bar z}^*$; no
orthonormality of $J$ is required. Replacing the anchor $z_*$ by a real
$z_\#\in\rho(A_0)$ replaces $M_z$ by $M_z-M_{z_\#}$ and $\Theta$ by
$\Theta-\Pi M_{z_\#}\Pi$ on $\ran\Pi$.
For an invertible change $\tau'=T\tau$, let $E=\ran\Pi$,
$E'=(T^*)^{-1}E$, and let $P_{E'}$ be the orthogonal projection
onto $E'$. The same extension has coordinates
\[
 \gamma'_z=\gamma_zT^*,\qquad M'_z=TM_zT^*,\qquad
 \Pi'=P_{E'},\qquad \Theta'=P_{E'}T\Theta T^*|_{E'} .
\]
In the full active space this gives $\Theta'=T\Theta T^*$.
\end{enumerate}
\end{theorem}
\begin{proof}
The adjoint domain, extension classification and resolvent in (i) are
Posilicano's finite-trace construction
\cite[Theorem~2.1, Corollary~3.2 and Theorems~4.2--4.3]{Pos}.
With the real anchor $w=z_*$ in \cite[(2.6)]{Pos}, the translation is
\[
 R^{\rm P}(z)=-\mathcal R_0(z),\quad G^{\rm P}(z)=-\gamma_z,\quad
 \Gamma^{\rm P}(z)=M_z,\quad \Theta^{\rm P}=-\Theta .
\]
It takes the negative of \cite[(2.7)]{Pos} to \eqref{fd:eq:projectionresolvent}.
Equations (2.4)--(2.5) of \cite{Pos} give (ii), since
$\gamma_z=(\tau\mathcal R_0(\bar z))^*$ is injective by surjectivity.

For (iii), write $J=UC$ with $U$ isometric and $C=(J^*J)^{1/2}>0$.
Part (i) applies on $\ran J$ with parameter $UC^{-1}\Theta_JC^{-1}U^*$.
Writing the singular coefficient as $Jc$, the boundary
condition is $(\Theta_J-J^*M_zJ)c=J^*\gamma_{\bar z}^*f$, which gives
the asserted correction and is invertible for nonreal $z$ by (ii).
At a new anchor $\Gamma_1$ becomes $\Gamma_1-M_{z_\#}\Gamma_0$.
Under $\tau'=T\tau$ the two boundary values $(c,d)$ become
$((T^*)^{-1}c,Td)$. Since $T(E^\perp)=(E')^\perp$, projecting
$d=\Theta c+w$, $w\in E^\perp$, gives exactly the stated $\Theta'$.
\end{proof}

\subsection{Point restrictions for a fixed outer realization}\label{fd:sec:found-outer}

\label{fd:sec:bounded}
\begin{theorem}[Complete interior jets for a regular elliptic realization]
\label{fd:thm:boundedjets}
Let $\Omega\subset\R^n$ be bounded and smooth, and let $A_0$ be a
self-adjoint realization of a formally symmetric elliptic differential
expression $P$ of order $2m$, with smooth coefficients. Assume
$C_c^\infty(\Omega)\subset\Dom(A_0)\subset H^{2m}(\Omega)$,
$A_0u=Pu$ there, and
\begin{equation}\label{fd:eq:graphregularity}
 c\norm u_{H^{2m}(\Omega)}\le
 \norm u_2+\norm{A_0u}_2\le C\norm u_{H^{2m}(\Omega)}.
\end{equation}
Let $X=\{x_1,\ldots,x_N\}\Subset\Omega$ be distinct.
Define $S_X$ by closing the restriction of $A_0$ to its domain elements
which vanish in a neighborhood of $X$. Put
$t=2m-n/2$ and $k=\lceil t\rceil-1$ when $t>0$, with an empty list
when $t\le0$. Then
\begin{equation}\label{fd:eq:boundedminimal}
 \Dom(S_X)=\{u\in\Dom(A_0):
       \partial^\alpha u(x_j)=0,\ |\alpha|\le k,\ 1\le j\le N\}.
\end{equation}
The complete jet map $\tau_X$ is graph-continuous and onto, and
$\ker\tau_X$ is dense in $L^2(\Omega)$. The deficiency indices are
$(Nr,Nr)$, with $r$ as in \eqref{fd:eq:count}. For every
$z\in\rho(A_0)$ its defect space is
\begin{equation}\label{fd:eq:boundeddefect}
 \ker(S_X^*-z)=\ran\gamma_X(z),\qquad
 \gamma_X(z)=(\tau_X(A_0-\bar z)^{-1})^*.
\end{equation}
The outer realization stays fixed in all these assertions.
In particular if $n\ge4m$, this restriction is just $A_0$.
\end{theorem}
\begin{proof}
By \eqref{fd:eq:graphregularity}, local cutoffs reduce the
annihilator calculation to Theorem~\ref{fd:thm:defect} in disjoint
interior balls. This gives \eqref{fd:eq:boundedminimal} with the
original outer conditions. Cutoff polynomials realize arbitrary jets,
and tests supported off $X$ give density of the kernel. Compact Sobolev
embedding gives compact resolvent and hence a real anchor.
Theorem~\ref{fd:thm:trace} now yields \eqref{fd:eq:boundeddefect}
and the index.
\end{proof}

\begin{corollary}[Concrete smooth bounded-domain realizations]
\label{fd:cor:realizations}
For $P=(-\Delta)^m$ on a bounded smooth domain, the theorem applies to
both of the following references:
\[
 \begin{array}{ll}
 \text{clamped:}&\Dom(A_0)=H^{2m}(\Omega)\cap H_0^m(\Omega),\\
 \text{Navier:}&A_0=(-\Delta_D)^m,\quad
 \Dom(A_0)=\{u\in H^{2m}(\Omega):
                   \Delta^j u|_{\partial\Omega}=0,\ 0\le j<m\}.
 \end{array}
\]
Both references are positive and have compact resolvent. Thus the
finite-index bounded-domain statements apply directly to these specified
reference operators.
\end{corollary}
\begin{proof}
The clamped domain and the graph estimate follow from
\cite[Corollary~2.21]{GGS} with $p=2$, $k=2m$. Iterating the
Dirichlet estimate for $-\Delta_D$ gives the stated Navier domain.
Positivity follows from the clamped form and the positive Dirichlet
Laplacian, respectively, and compactness from Sobolev embedding.
\end{proof}

\subsection{Lower bounds and eigenvalue counts}\label{fd:sec:found-spectrum}

\label{fd:sec:spectral}
\begin{theorem}[Finite index bounds the number of lower eigenvalues]
\label{fd:thm:spectral}
In the finite-trace framework let $A_0\ge aI$ and let the trace
dimension be $q$. Every self-adjoint extension $A_L$ of $S$ is bounded
below, and its spectral subspace below $a$ has dimension at most $q$.
If $A_0$ has compact resolvent, so does $A_L$, and
\begin{equation}\label{fd:eq:counting}
 |N_L(t)-N_0(t)|\le q,\qquad
 N_j(t)=\dim\operatorname{ran}\boldsymbol1_{(-\infty,t)}(A_j).
\end{equation}
Under this compact-resolvent assumption, both counting functions are
finite for every real $t$.
\end{theorem}
\begin{proof}
The finite-deficiency gap bound \cite[Section~8.3]{Wei} gives at most
$q$ dimensions in any bounded open subinterval of $(-\infty,a)$.
Exhaustion gives the same bound below $a$, hence semiboundedness.
Compactness passes through the finite-rank resolvent formula
\eqref{fd:eq:projectionresolvent}. In the compact case intersect the
spectral subspace below $t$ of $A_L$ with $\ker\Gamma_0$, or that of
$A_0$ with $\ker\tau$. Each intersection loses at most $q$ dimensions,
and the two operators agree on it. The variational dimension bound
\cite[Theorem~4.15(i)]{Teschl} gives both inequalities in
\eqref{fd:eq:counting}.
\end{proof}

Uniformity across extensions requires a separate argument. For example,
$b=\alpha^{-1}$ in \eqref{eq:res} gives a one-site eigenvalue
$-\alpha^2/64$ for $\alpha<0$, which escapes to $-\infty$ if $\alpha=-1/h$.

\section{The repaired domain and the planar Weyl matrix}\label{fd:sec:found-repairmatrix}

\subsection{The repaired singular domain and its \texorpdfstring{$L^2$}{L2} part}\label{fd:sec:found-repair}

\label{fd:sec:repair}
Return temporarily to $P=\Delta^m$ and its zero-energy Navier
realization $A_0=(-1)^m(-\Delta_D)^m$ on a bounded smooth domain.
Zero is in $\rho(A_0)$. Let $T_1,\ldots,T_d$ be independent
point-supported distributions with derivative degrees at most $2m-1$,
and let $\psi_j$ be their zero-energy Navier Green fields.
For the family of the preceding flux section these are
$\Delta_y^jG$ and $\partial_{y_s}\Delta_y^jG$.

\begin{theorem}[Direct-sum repair and the exact $L^2$ criterion]
\label{fd:thm:repair}
The fields $\psi_j$ have locally integrable representatives and
$P\psi_j=T_j$ distributionally. The space
\begin{equation}\label{fd:eq:repairedambient}
 \mathcal E_T=\left\{u_0+\sum_{j=1}^dc_j\psi_j:
                  u_0\in\Dom(A_0),\ c\in\C^d\right\}
\end{equation}
is a direct sum. Equip it with the squared norm
$\norm{u_0}_{\Dom(A_0)}^2+|c|^2$ and define the singular coordinate
$\beta(u)=c$. The punctured differential expression on it is
$P_{\rm off}u=A_0u_0$.
For any $f\in L^2$ and $c\in\C^d$ the conditions
$P_{\rm off}u=f$, $\beta(u)=c$ have the unique solution
\begin{equation}\label{fd:eq:repairedsolution}
 u=A_0^{-1}f+\sum_jc_j\psi_j.
\end{equation}
It depends boundedly on $(f,c)$ in the declared norm.

At each puncture combine $\sum_jc_jT_j$ into a derivative polynomial
and let $E_{\rm ad}\subset\C^d$ be the subspace whose degree at every
point is strictly less than $2m-n/2$. Then
\begin{equation}\label{fd:eq:ambientL2part}
 \mathcal E_T\cap L^2
       =\Dom(A_0)\dotplus\left\{\sum_jc_j\psi_j:c\in E_{\rm ad}\right\}.
\end{equation}
In particular, replacing $c$ in \eqref{fd:eq:repairedsolution} by any
bounded functional vector $Cf$ gives an $L^2$-valued solution map for
every $f$ if and only if $\ran C\subset E_{\rm ad}$.
\end{theorem}
\begin{proof}
For the constant-coefficient operator choose a fundamental solution
$\Psi$ and an interior cutoff in $x-y$ equal to one near zero.
Applying $P$ to the cut-off fundamental solution gives $\delta_y$
plus a smooth remainder supported away from $y$. Subtract the Navier
inverse of that remainder. Differentiating in $y$ constructs every
$\psi_j$ with the prescribed distributional source and outer conditions.
For $0\le j\le 2m-1$, the order-$j$ derivatives of the
fundamental solution off the pole are
$O(r^{2m-n-j}(1+|\log r|))$, hence locally integrable.
In identifying the $j$th distributional derivative by integration
by parts for $1\le j\le2m-1$, the inner boundary term from an
order-$(j-1)$ derivative is $O(r^{2m-j}(1+|\log r|))$ for each
fixed test function. It tends to zero. Induction identifies the
distributional derivatives with these locally integrable fields.

If a combination $\sum_jc_j\psi_j$ belongs to $\Dom(A_0)$, its
distributional $P$ image is both $L^2$ and the point-supported
distribution $\sum_jc_jT_j$. Such a distribution is zero, and
independence of the $T_j$ gives $c=0$. This proves directness.
Applying $P$ off the punctures then proves the solution formula and
uniqueness. Completeness and boundedness follow from the product norm.

At each pole, the local test in Proposition~\ref{fd:prop:channels}
applies to the combined source $\sum_jc_jT_j$, including cancellations
of leading derivatives. Away from the poles the fields are smooth.
This proves \eqref{fd:eq:ambientL2part} and the condition on $\ran C$.
\end{proof}

For the radial--gradient list, the normalized functionals in
Proposition~\ref{fd:prop:flux} recover precisely the coordinates $c_j$
on \eqref{fd:eq:repairedambient}. They vanish on regular $H^{2m}$
summands, not only on smooth summands: the divergence theorem writes
their flux as an integral of $\Delta^{m-i}u_0$ over a shrinking ball,
and their first moment as an integral of
$\partial_s\Delta^{m-i-1}u_0$. Both tend to zero by Cauchy--Schwarz.
The spherical integrals have their usual Sobolev trace meaning for
positive radius.

Thus a correct zero-energy version of the pattern of
\cite[Theorems~2.1--2.2]{KT} is available in a specified finite
singular direct sum. In $L^2$ one either chooses an admissible sublist
as an explicit submodel or uses the complete jet list of
Theorem~\ref{fd:thm:boundedjets}. In the planar biharmonic case the
complete list has six channels, with two trace-free Hessians replacing
the two inadmissible derivatives of $\Delta G$.
An $L^2$ resolvent requires $L^2$ source fields and a dense operator
domain. For self-adjoint relations these requirements are supplied by
Theorem~\ref{fd:thm:trace}; no continuation of a
reference resolvent to an unspecified distribution space is used.

\subsection{The complete planar biharmonic Weyl matrix}\label{fd:sec:found-weyl}

\label{fd:sec:planar}
We make the complete six-channel case explicit. Take
$A_0=\Delta^2$ on $H^4(\R^2)$, the anchor $z_*=-1$, and
\[
 \tau u=(u(0),\partial_1u(0),\partial_2u(0),
                   \partial_1^2u(0),\partial_1\partial_2u(0),
                   \partial_2^2u(0))^T.
\]
For $z=-\lambda$ set
$\widehat G_\lambda=(|\xi|^4+\lambda)^{-1}$ and
\[
 \gamma_{-\lambda}=(G_\lambda,-\partial_1G_\lambda,
 -\partial_2G_\lambda,\partial_1^2G_\lambda,
 \partial_1\partial_2G_\lambda,\partial_2^2G_\lambda).
\]
These signs are the adjoints of the stipulated regular jets.
Write
\[
 g_\lambda=\frac1{8\sqrt\lambda},\qquad
 L_\lambda=\frac{\log\lambda}{16\pi},\qquad
 T_\lambda=\frac{1-\sqrt\lambda}{64}.
\]
Put $U_D=(G_\lambda,-\partial_eG_\lambda)$ and
$U_E=(G_\lambda,\partial_e^2G_\lambda)$.
The square root and logarithm have their principal branches on
$\C\setminus(-\infty,0]$.

\begin{theorem}[Complete six-channel matrix and two active collision families]
\label{fd:thm:planarM}
For the stated regular jet basis,
\begin{equation}\label{fd:eq:fullM}
 M_{-\lambda}=\begin{pmatrix}
 g_\lambda-g_1&0&0&L_\lambda&0&L_\lambda\\
 0&-L_\lambda&0&0&0&0\\
 0&0&-L_\lambda&0&0&0\\
 L_\lambda&0&0&3T_\lambda&0&T_\lambda\\
 0&0&0&0&T_\lambda&0\\
 L_\lambda&0&0&T_\lambda&0&3T_\lambda
 \end{pmatrix}.
\end{equation}
\par\Needspace{9\baselineskip}
For any real unit direction $e$, the smaller traces
$\tau_Du=(u(0),\partial_eu(0))$ and
$\tau_Eu=(u(0),\partial_e^2u(0))$ have Weyl matrices
\begin{align}
 M_D(-\lambda)&=\begin{pmatrix}
 g_\lambda-g_1&0\\0&-L_\lambda
 \end{pmatrix},\label{fd:eq:MD}\\
 M_E(-\lambda)&=\begin{pmatrix}
 g_\lambda-g_1&L_\lambda\\L_\lambda&3T_\lambda
 \end{pmatrix}.\label{fd:eq:ME}
\end{align}
For real $a,\theta,\beta,\vartheta$, the Hermitian active parameters
\begin{equation}\label{fd:eq:activeTheta}
 \Theta_D=\begin{pmatrix}-g_1&a\\a&-\theta\end{pmatrix},\qquad
 \Theta_E=\begin{pmatrix}-g_1&2\beta\\2\beta&-\vartheta\end{pmatrix}
\end{equation}
therefore define self-adjoint extensions of the complete six-jet
restriction, whose resolvents are respectively
\begin{align}
 R_D(-\lambda)&=\mathcal R_0(-\lambda)
  -U_D
  \begin{pmatrix}g_\lambda&-a\\-a&\theta-L_\lambda\end{pmatrix}^{-1}
  U_D^*,\label{fd:eq:Dfamily}\\
 R_E(-\lambda)&=\mathcal R_0(-\lambda)
  -U_E
  \begin{pmatrix}g_\lambda&-2\beta+L_\lambda\\
             -2\beta+L_\lambda&\vartheta+3T_\lambda\end{pmatrix}^{-1}
  U_E^* .\label{fd:eq:Efamily}
\end{align}
The last two displays are written for real $\lambda>0$ whenever their
matrices are invertible. For complex $\lambda$ the adjoint source
factor has parameter $\bar\lambda$.
\end{theorem}
\begin{proof}
With $F_\lambda=(|\xi|^4+\lambda)^{-1}-(|\xi|^4+1)^{-1}$, the multipliers
$(-i\xi)^\beta F_\lambda$, $|\beta|\le2$, are $O(|\xi|^{|\beta|-8})$,
so the source differences lie in $H^4$ and the entries
$M_{\alpha\beta}=\partial_x^\alpha(-\partial_x)^\beta(G_\lambda-G_1)(0)$
are ordinary traces, although separate Hessian traces of $G_\lambda$
diverge. Thus $M_{\alpha\beta}$ is $i^{|\alpha|}(-i)^{|\beta|}$ times
the moment $(2\pi)^{-2}\int\xi^{\alpha+\beta}F_\lambda\,\dd\xi$. Odd
moments vanish, the zeroth is $g_\lambda-g_1$, the quadratic and pure
quartic moments are $I_2(\lambda)=-L_\lambda$ and $I_4(\lambda)=3T_\lambda$
by \eqref{eq:I2}--\eqref{eq:I4}, and the mixed quartic moment is
$I_4/3$, since $\int_0^{2\pi}\cos^2\theta\sin^2\theta\,\dd\theta
=\frac13\int_0^{2\pi}\cos^4\theta\,\dd\theta$. This gives
\eqref{fd:eq:fullM}. Compression with the linear and quadratic tensors of
$e$ (the mixed coefficient being $2e_1e_2$) gives
\eqref{fd:eq:MD}--\eqref{fd:eq:ME}. Finally $\Theta_D-M_D$ and
$\Theta_E-M_E$ are the negatives of the matrices inverted in
\eqref{fd:eq:Dfamily}--\eqref{fd:eq:Efamily}, so
Theorem~\ref{fd:thm:trace}(iii) gives these resolvents; the
complete six-jet kernel lies in each smaller-trace kernel.
\end{proof}

For a bounded Navier polyharmonic reference, a zero-energy local Green
kernel is $\Psi(x-y)+H(x,y)$ with $H$ smooth jointly near an interior
pole. The cutoff construction in Theorem~\ref{fd:thm:repair} proves this:
its smooth remainder depends smoothly on $y$, and smooth elliptic
regularity for its Navier inverse gives joint interior derivatives of
every order. The complete local jet list is consequently the same as
in the whole space, but its regularized matrices depend on the outer
realization. The abstract difference definition
$M_z=\tau(\gamma_z-\gamma_*)$ supplies finite coordinates directly,
without evaluating separate divergent point traces. Local agreement
of singularities does not by itself prove a uniform collision
approximation or its common lower bound in the bounded domain.

\section*{Declarations}

\subsection*{Funding}
This research was funded by the Committee of Science of the Ministry of
Science and Higher Education of the Republic of Kazakhstan (Grant
Nos.~BR31714735 and AP32319153).

\subsection*{Conflict of interest}
The author declares that there is no conflict of interest.

\subsection*{Data availability}
No datasets were generated or analysed during the present study.



\subsection*{Declaration of generative AI use}
The author used ChatGPT (OpenAI) and Claude (Anthropic) for language
editing. ChatGPT also assisted with literature searches, verification of
references, revision of proofs, and symbolic and numerical checks of
selected calculations. Claude checked the proofs and proposed corrections
and revisions of several arguments, including Lemma~\ref{lem:pseudo},
Remark~\ref{rem:slowonset}, the statement on the centre tuning in
Theorem~\ref{thm:twohessians} and the statement on the pinned monopole in
Theorem~\ref{thm:obtusedipoles}. Claude also carried out symbolic and
numerical checks and verified bibliographic data. The author reviewed and
edited all such output and takes full responsibility for the proofs,
citations and final text.

\end{document}